\documentclass[10pt,letterpaper]{article}
\usepackage[top=0.85in,left=2.75in,footskip=0.75in,marginparwidth=2in]{geometry}

\usepackage[utf8]{inputenc}
\usepackage{cite}
\usepackage{nameref,hyperref}

\usepackage{microtype}
\DisableLigatures[f]{encoding = *, family = * }

\raggedright
\usepackage{changepage}
\usepackage[aboveskip=1pt,labelfont=bf,labelsep=period,singlelinecheck=off]{caption}

\makeatletter
\renewcommand{\@biblabel}[1]{\quad#1.}
\makeatother

\usepackage{lastpage,fancyhdr,graphicx}
\usepackage{epstopdf}
\fancyheadoffset[L]{2.25in}
\fancyfootoffset[L]{2.25in}

\usepackage{color}
\definecolor{Gray}{gray}{.25}

\usepackage{graphicx}
\usepackage{sidecap}
\usepackage{wrapfig}
\usepackage[pscoord]{eso-pic}
\usepackage[fulladjust]{marginnote}
\reversemarginpar

\usepackage{amsmath,amssymb,amsfonts,amstext,mathtools,bm,latexsym}
\usepackage{siunitx}
\usepackage{float}
\usepackage{booktabs}
\usepackage{multirow}
\usepackage{subcaption}
\usepackage{algorithm}
\usepackage{algpseudocode}
\usepackage{todonotes}

\newcommand{\ignore}[1]{}  
\newcommand{\vect}[1]{\bm{#1}}
\newcommand{\setS}{\mathcal{S}}

\newcommand{\Cstar}{\mathcal{C}^{\ast}}
\newcommand{\KL}{\mathcal{K}_\infty}

\DeclareMathOperator*{\argmin}{arg\,min}

\begin{document}
\vspace*{0.35in}

\begin{flushleft}
{\Large
\textbf\newline{Learning Safety-Guaranteed, Non-Greedy Control Barrier Functions Using Reinforcement Learning}
}
\newline
\\
Minduli Wijayatunga\textsuperscript{1,*},
Nathan Wallace\textsuperscript{2},
Salah Sukkarieh\textsuperscript{2},
Roberto Armellin\textsuperscript{3}
\\
\bigskip
\bf{1} University of Illinois Urbana-Champaign, Urbana, IL 61801, USA
\\
\bf{2} Australian Centre for Robotics, The University of Sydney, NSW 2006, Australia
\\
\bf{3} University of Auckland, Auckland 1010, New Zealand
\\
\bigskip
* corresponding: minduli@illinois.edu

\end{flushleft}

\section*{Abstract}
Spacecraft rendezvous and proximity operations (RPO) present inherent safety risks to high-value assets, making safety guarantees essential for mission success. However, overly conservative control policies developed for safety can reduce mission efficiency.  This work proposes a unified two-stage reinforcement learning (RL) framework that addresses the complementary limitations of a traditional safety approach—Input Constrained Control Barrier Functions (ICCBFs)—in safety-critical, fuel-limited spacecraft control. For a certified safe set $\setS$, ICCBFs yield a provably invariant inner set $\Cstar\subseteq \setS$ under input bounds, but the resulting per-step Quadratic Program (QP) is greedy and fuel-inefficient inside $\Cstar$, and states in the residual set $\setS\setminus\Cstar$ that are recoverable are conservatively abandoned. {Stage 1} of the proposed framework learns state-dependent class-$\KL$ parameters that adapt the ICCBF/CLF decay rates to embed long-horizon cost awareness while maintaining invariance in $\Cstar$. {Stage 2} then learns a residual barrier $h_{\mathrm{RL}}(\vect{x})$ that provides recoverability for a subset of $\setS\setminus\Cstar$. At run-time, the framework selects the appropriate barrier formulation (Stage 1 or Stage 2) and solves a lightweight QP under zero-order hold.  Both stages are trained with PPO using rewards that penalise constraint violations, control effort, and task-specific metrics.  The framework is evaluated on three problems: cruise control, spacecraft rendezvous with a rotating target, and a spacecraft inspection task that maximizes observability while respecting keep-in and keep-out zone constraints. 
The test cases show a reduction in median fuel consumption compared to the ICCBF baselines by \SIrange{25}{12}{\percent}, and a \SIrange{7}{8}{\percent} increase in the number of cases that remain inside the safe set $\setS$. These results show that RL can embed long-term cost awareness into safety-critical control structures while maintaining the computational efficiency of QPs.


\section{Introduction}
Autonomous spacecraft proximity operations such as docking, inspection, and formation flying often operate in complex and dynamic environments where collision avoidance, line-of-sight maintenance, approach corridor, and keep-out zone (KOZ) adherence are critical for mission success and asset protection. As such, they demand rigorous safety guarantees while operating under constraints in fuel consumption, computational resources, and real-time response requirements. Traditional guidance and control approaches often struggle to balance the competing demands of safety and operational constraints, particularly in scenarios involving simultaneous constraints and multiple objectives. The challenge is compounded by the need for real-time decision-making capabilities that can adapt to changing mission requirements while providing and adhering to mathematical guarantees of safety. This combination has proven difficult to achieve with conventional control methods.

For a given system, safety is usually certified by defining a safe subset in its state space, within which the system states must remain \cite{agrawal2021safe}. Control Barrier Functions (CBFs) have emerged as a promising set-theoretic approach for addressing safety assurance, offering a mathematically rigorous framework for safety-critical control in control-affine systems \cite{xiao2023safe}. CBFs provide several key advantages in safety enforcement:  they enable real-time computation through quadratic programming (QP) formulations and allow compositional safety by combining multiple safety constraints. Implementation of CBFs via QP introduces a safety filter that modifies desired control inputs to ensure forward invariance of safe sets, making it particularly attractive for spacecraft applications where safety cannot be compromised. Thus, CBF--QP's computational efficiency and theoretical safety guarantees have led to its successful application in various domains in robotics and aerospace \cite{8796030}. 


However, the CBF-QP formulation alone cannot provide overall fuel optimality or adherence to performance criteria beyond safety. Traditional CBF-QP formulations suffer from inherently greedy behavior that prioritizes immediate constraint satisfaction without considering long-term trajectory optimality or fuel efficiency. The standard QP problem formulation acts greedily at each timestep, often leading to oscillatory control behavior, excessive thrust consumption over time, and suboptimal trajectories in spacecraft proximity operations. This greediness becomes particularly problematic in time-varying scenarios such as docking with rotating targets, where the instantaneous safety-preserving control may conflict with the globally optimal approach strategy.  A further complication arises when the original CBF expression lacks velocity dependence, causing the Lie derivative to approach zero and preventing effective safety enforcement. These limitations highlight a fundamental trade-off: while CBFs excel at providing safety guarantees, they do not inherently optimize for mission-critical objectives such as fuel minimization, trajectory smoothness, or task-specific performance metrics.

Recent research has attempted to address these limitations. Enhanced barrier constructions such as Higher-Order CBFs (HOCBFs) \cite{8796030} and  ICCBFs \cite{agrawal2021safe} have addressed the velocity independence issue by incorporating higher-order derivatives or input constraints, though they do not resolve the fundamental greediness problem. ICCBFs introduce additional constraints that limit the feasible domain, potentially leading to infeasible barrier constructions and subsequent QP infeasibility that causes control system failure. 
Some formulations embed CBF constraints within Model Predictive Control (MPC) frameworks \cite{zeng2021enhancing,9483029}, or employ predictive safety filters to reduce greedy behavior of CBFs \cite{breeden2022predictive}. These methods can offer reduced greediness at the expense of higher computational load and complex tuning, which can be problematic at the high control rates present in spacecraft problems. Recently, data-driven approaches that parametrize barrier functions using neural networks \cite{10015199} have also been explored, where the parameters are learned from demonstrations or reinforcement learning (RL). These methods either constrain the neural network’s actions using a pre-defined safety certificate \cite{cheng2019end,berkenkamp2017safe},  learn safety certificates alongside a control policy \cite{westenbroek2021combining}, or derive certificates directly from observable quantities \cite{chow2018lyapunov}. While effective, these works primarily focus on constructing valid CBFs through learning, rather than on learning CBFs that additionally account for the optimization of other objectives, such as fuel usage. 

In contrast, this paper focuses on developing non-greedy CBFs through a unified RL-ICCBF framework that preserves the lightweight, certifiable CBF-QP structure at runtime. The approach addresses ICCBF's limitations through a coordinated two-stage framework: Stage 1 (RL-tunned ICCBF) learns slack parameters applied to existing ICCBF constraints to enable non-greedy control within the proven-safe inner set $C^*$ (defined in Section \ref{pformcbf}), while Stage 2 (RL-residual CBF) learns a residual barrier function $h_{RL}(x)$ that recovers some states in the residual set $\mathcal{R} = \mathcal{S} \setminus C^*$ that ICCBF conservatively abandons. In this stage, the gradients of the learned barrier function are obtained via automatic differentiation. Both stages use Proximal Policy Optimization (PPO) \cite{schulman2017proximal} to minimize cumulative fuel consumption, other optional objectives, and safety violations over complete trajectories, thereby embedding long-horizon cost awareness while maintaining per-step safety guarantees. The framework's effectiveness is compared against ICCBF by running the benchmarks presented in the ICCBF paper \cite{agrawal2021safe}, which includes cruise control and spacecraft rendezvous with a rotating target. Then, the framework is applied to an additional three-dimensional inspection task where an observability/lighting metric is optimised while adhering to keep-in/keep-out constraints and thrust limits.

In summary, the main contributions of this work are:
\begin{itemize}    
    \item An RL-tuned ICCBF approach that learns state-dependent class-$\mathcal{K}_\infty$ functions that enable non-greedy control within the proven-safe inner set $C^*$.
    
    \item An RL-residual CBF method that recovers states in the residual set $\mathcal{R} = \mathcal{S} \setminus C^*$ that ICCBF conservatively abandons.
    
    \item A comprehensive framework that maintains the real-time computational advantages of QP while reducing greedy behaviour and expanding the feasible set through RL, making it applicable for safety-critical spacecraft applications such as autonomous docking and inspection,
\end{itemize}

The remainder of this paper is structured as follows. Section \ref{pformcbf} presents the theoretical background on 
CBFs, ICCBF and RL. Section \ref{methodology} details the methodology of the unified RL framework.  Section \ref{results_benchmark} discusses the implementation of the developed framework on two benchmark problems from the literature, followed by Section \ref{application}, which extends the framework to optimize the inspection of an uncooperative space object while maintaining safety in three-dimensional space.  Finally, Section \ref{conclusions} concludes the paper by analysing the performance of the framework over the three test cases and discussing avenues of future work.

\section{Theoretical Background} \label{pformcbf}

\textbf{Notation:} In this work, the set of real numbers is denoted $\mathbb{R}$ and non-negative real numbers are denoted as $\mathbb{R}_+$. A continuous function $\alpha: [0, a) \to [0, \infty)$ is a class-$\mathcal{K}_{\infty}$ function if it is strictly increasing and $\alpha(0) = 0$. The Lie derivative of $h(x)$ along $f(x)$ is denoted $L_f h(x) = \frac{\partial h}{\partial x} f(x)$. Similarly, $L_g h(x) = \frac{\partial h}{\partial x} g(x)$. Note that Int($\mathcal{C}$) denotes the interior of a set $\mathcal{C}$, while $\delta \mathcal{C}$ denotes the boundary of it. Note that the notation and terminology in this work follow \cite{agrawal2021safe} to maintain consistency with the foundational ICCBF framework.

\subsection{Problem Setup}
Consider a control-affine nonlinear dynamical system with state $x \in \mathcal{X} \subset \mathbb{R}^n$ and control input $u \in \mathcal{U} \subset \mathbb{R}^m$:
\begin{equation}\label{xdoteq}
\dot{x} = f(x) + g(x)u,
\end{equation}
where $f : \mathcal{X} \to \mathbb{R}^n$ and $g : \mathcal{X} \to \mathbb{R}^{n \times m}$ are sufficiently smooth functions.  Under a Lipschitz continuous feedback law $u = \pi(x)$, the closed-loop system becomes:
\begin{equation}
\dot{x} = f(x) + g(x)\pi(x)
\end{equation}

A state $x$ is defined as safe if it lies in the set $\mathcal{S}$, the zero superlevel set of the function  $h : \mathcal{X} \to \mathbb{R}$, which is continuously differentiable. 
\begin{align}
\mathcal{S} &= \{x \in \mathcal{X} : h(x) \geq 0\} \\
\partial\mathcal{S} &= \{x \in \mathcal{X} : h(x) = 0\} \\
\text{Int}(\mathcal{S}) &= \{x \in \mathcal{X} : h(x) > 0\}
\end{align}
Here, $\mathcal{S}$ is the safe set, assumed to be closed, non-empty, and simply connected.

\textbf{Definition 1 (Forward Invariance).} A set $\mathcal{S}$ is rendered forward invariant by a feedback controller $\pi : \mathcal{S} \to \mathcal{U}$ if, for the closed-loop system, $x(0) \in \mathcal{S}$ implies $x(t) \in \mathcal{S}$ for all $t \geq 0$.

\subsection{Control Barrier Functions (CBFs)}

Traditional CBF theory provides a framework for ensuring forward invariance through real-time optimization. The fundamental result connecting barrier functions to safety is given by Nagumo's theorem:

\textbf{Lemma 1 (Nagumo's Theorem).} Consider the system with safe set $\mathcal{S}$ defined by a continuously differentiable function $h : \mathcal{X} \to \mathbb{R}$. For a Lipschitz continuous feedback controller $\pi : \mathcal{S} \to \mathcal{U}$, the set $\mathcal{S}$ is forward invariant if and only if:
\begin{equation}
L_f h(x) + L_g h(x) \pi(x) \geq 0, \quad \forall x \in \partial\mathcal{S}
\end{equation}

\textbf{Definition 2 (CBF).} Let $\mathcal{S} \subset \mathcal{X} \subset \mathbb{R}^n$ be the superlevel set of a continuously differentiable function $h : \mathcal{X} \to \mathbb{R}$. The function $h$ is a CBF if there exists an extended class-$\mathcal{K}_\infty$ function $\alpha$ such that:
\begin{equation}
\sup_{u \in \mathcal{U}} [L_f h(x) + L_g h(x) u] \geq -\alpha(h(x)), \quad \forall x \in \mathcal{X}
\end{equation}

Note that this is the definition provided in \cite{agrawal2021safe}. A stronger notion of CBFs are introduced in  \cite{8796030} and \cite{7782377}.

\textbf{Theorem 1 (CBF Safety).} If $h$ is a CBF and $\frac{\partial h}{\partial x}(x) \neq 0$ for all $x \in \partial\mathcal{S}$, then any Lipschitz continuous controller $\pi(x) \in K_{\text{CBF}}(x)$, where:
\begin{equation}
K_{\text{CBF}}(x) = \{u \in \mathcal{U} : L_f h(x) + L_g h(x) u + \alpha(h(x)) \geq 0\}
\end{equation}
renders the set $\mathcal{S}$ forward invariant.

\subsection{Control Lyapunov Functions (CLF)}
Performance objectives in control systems are often expressed through the concept of asymptotic stability. A CLF provides a constructive approach to ensure convergence to a desired equilibrium or set. CLFs are commonly used in conjunction with CBFs to stabilize a system to an equilibrium point while maintaining safety. 

\textbf{Definition 3 (CLF).} A continuously differentiable function $V : \mathcal{X} \to \mathbb{R}$ is a CLF if there exists a class-$\mathcal{K}_\infty$ function $\alpha'$ such that for all $x \in \mathcal{X}$:
\begin{equation}\label{eq:clf_condition}
    \inf_{u \in \mathcal{U}} [L_f V(x) + L_g V(x) u] \leq -\alpha'(V(x)) 
\end{equation}
The existence of a CLF yields a family of controllers that asymptotically stabilize the system. The set of control inputs satisfying the CLF condition is:
\begin{equation}
K_{\text{CLF}}(x) = \{u \in \mathcal{U} : L_f V(x) + L_g V(x) u + \alpha'(V(x)) \leq 0\}
\end{equation}

Any locally Lipschitz controller $u(x) \in K_{\text{CLF}}(x)$ renders the origin asymptotically stable for the closed-loop system.

\subsection{CLF-CBF-QP Formulation}

In practice, CBF constraints are often combined with CLF objectives in a QP framework. Given a CLF $V(x)$ that encodes convergence to a desired target, the standard CLF-CBF-QP formulation is:
\begin{align}
u^* = \argmin_{u \in \mathbb{R}^m, \delta \in \mathbb{R}_+} \quad & \frac{1}{2} u^T u + p_1 \delta^2 \label{eq:qp_obj}\\
\text{subject to} \quad 
& L_f h(x) + L_g h(x) u \geq -\alpha(h(x)) \label{eq:cbf_constraint}\\
& L_f V(x) + L_g V(x) u \leq -\alpha' (V(x))+ \delta \label{eq:clf_constraint}\\
& u \in \mathcal{U} \label{eq:input_constraint}
\end{align}
where $p_1 > 0$ is a penalty weight, and $\delta$ provides relaxation for CLF feasibility. $\alpha$ and $\alpha'$ are class-$\mathcal{K}_\infty$ functions. Note that in this paper, the test cases presented in Section \ref{results_benchmark} utilizes CLF--CBF--QP, while the case in Section \ref{application} uses a CBF--QP formulation.

\subsection{Limitations of CBF}
The traditional CBF implementations can exhibit the following limitations. 

\begin{enumerate}
    \item \emph{When $h(x)$ has relative degree greater than one with respect to the control input  (i.e., $L_g h(x) = 0$), the standard CBF condition becomes ineffective}. This issue arises frequently in aerospace applications where safety constraints may not directly depend on velocities. For example, consider a line-of-sight constraint in spacecraft docking:
    \begin{equation}
        h(x) = \frac{\mathbf{r}_{cp} \cdot \mathbf{e}}{\|\mathbf{r}_{cp}\|} - \cos(\gamma),
    \end{equation}
    where $\mathbf{r}_{cp}$ is the relative position vector and $\mathbf{e}$ is the docking axis. Since $h(x)$ depends only on position and orientation (not velocities), $L_g h(x) = 0$, rendering the CBF constraint ineffective.
    
    \item Standard CBF theory \emph{presumes sufficient control authority to render the safe set forward invariant}. However, input constraints $u \in \mathcal{U}$ may violate this assumption, leading to infeasible QP formulations, reduced feasibility domains, and the consequent loss of safety guarantees. 

        \item The standard CBF-CLF-QP formulation is \emph{inherently greedy}. It minimises instantaneous control effort while enforcing immediate safety constraints, without regard for long-term trajectory optimality. This can result in excessive fuel consumption, control chattering, and an inability to anticipate constraint violations that could otherwise be avoided with foresight.

\end{enumerate}

\subsection{Input Constrained Control Barrier Functions (ICCBF)}\label{ICCBFdetail}

When input constraints limit the controller authority, the safe set $\mathcal{S}=\{x:h(x)\ge 0\}$ obtained by a nominal safety function $h$ may not be forward invariant. The ICCBF construction generates an inner safe set $\mathcal{C}^\star\subseteq\mathcal{S}$ that is forward invariant with the available constrained inputs, by iteratively composing $h$ with Lie derivatives and class-$\mathcal{K}_\infty$ margins.

The ICCBF construction begins by recursively defining 
\begin{align}
b_0(x) &\triangleq h(x), \label{eq:iccbf-b0}\\
b_1(x) &\triangleq \inf_{u\in\mathcal{U}}\big[L_f b_0(x)+L_g b_0(x)\,u+\alpha_0\!\left(b_0(x)\right)\big], \label{eq:iccbf-b1}\\
&\ \ \vdots \nonumber\\
b_N(x) &\triangleq \inf_{u\in\mathcal{U}}\big[L_f b_{N-1}(x)+L_g b_{N-1}(x)\,u+\alpha_{N-1}\!\left(b_{N-1}(x)\right)\big],
\label{bees}
\end{align}
where each $\alpha_i$ is a class-$\mathcal{K}_\infty$ function and $N\in\mathbb{N}$ is chosen so that all expressions are well-defined and continuous on the domain of interest. For each level, the following superlevel sets are introduced 
\begin{equation}
\mathcal{C}_i \triangleq \{x : b_i(x)\ge 0\},\quad i=0,\dots,N,
\end{equation}
and the inner safe set is defined as the intersection 
\begin{equation}
    \mathcal{C}^\star \triangleq \bigcap_{i=0}^{N}\mathcal{C}_i.
    \label{cstar}
\end{equation}
Function $b_N$ is then an ICCBF on $\mathcal{C}^\star$ if there exists a class-$\mathcal{K}_\infty$ function $\alpha_N$ such that
\begin{equation}
\sup_{u\in\mathcal{U}}\big[L_f b_N(x)+L_g b_N(x)\,u+\alpha_N\!\left(b_N(x)\right)\big]\ge 0,\quad \forall x\in\mathcal{C}^\star.
\label{eq:iccbf-cond}
\end{equation}
Then any locally Lipschitz feedback $\pi(x)$ that satisfies the pointwise linear inequality
\begin{equation}
L_f b_N(x)+L_g b_N(x)\,\pi(x)+\alpha_N\!\left(b_N(x)\right)\ \ge\ 0,
\label{eq:iccbf-constraint}
\end{equation}
renders $\mathcal{C}^\star$ forward invariant (by Nagumo's theorem applied on each active boundary in the chain).

Limitations (1) and (2) can be addressed through ICCBFs. 
However, ICCBFs alone do not improve the greedy behavior in limitation (3) that is inherent in CBF-based control. To address this limitation, this work proposes a unified ICCBF-RL framework that learns non-greedy barrier functions, embedding long-term cost awareness directly into the safety certification process while preserving computational efficiency and safety guarantees.

\subsection{Reinforcement Learning (RL)} 

RL is a branch of machine learning where no reliance on a static dataset is needed. Instead, RL operates on a dynamic environment to learn by experience. 
The agent, policy, and reward/cost are critical to understanding how RL works. Given an environment, and the agent's state $x$ in it, a policy $\pi_\theta(a|x)$ determines what action $a$ RL will take. The RL agent explores the environment while rewards $r(x,a)$ provide feedback for a given action during the training process.  During training, an agent will take actions based on an untrained policy, and over time, this policy will improve until the optimal policy is learnt \cite{Sutton1998}.

For a sequence of states, the performance is evaluated by the discounted reward to go $\mathcal{R}^{\pi}_i$, which is  defined as 
\begin{equation}
\mathcal{R}^{\pi}_i = \sum_{j=0}^{\infty} \gamma^{j-1} , r^{\pi}(x_{j},a_{j}), \qquad \gamma\in(0,1],
\end{equation}
The agent's goal is to obtain a policy that maximises $\mathcal{R}^{\pi}_i$.

In this work, the RL algorithm PPO is used to develop and update CBFs, due to the continuous nature of the states and actions involved. PPO is known to be robust and effective in such environments, as it can efficiently navigate continuous action spaces by using a policy gradient approach while iteratively adjusting the policy based on sampled trajectories \cite{schulman2017proximal,10.1007/978-3-031-25755-1_9}. Furthermore, PPO allows clipping of the objective
function, which ensures a more stable training process \cite{HarryThesis}. 

\section{Methodology}\label{methodology}

\begin{figure*}[t]
  \centering
  \begin{subfigure}[t]{0.5\textwidth}
    \centering
    \includegraphics[width=\linewidth]{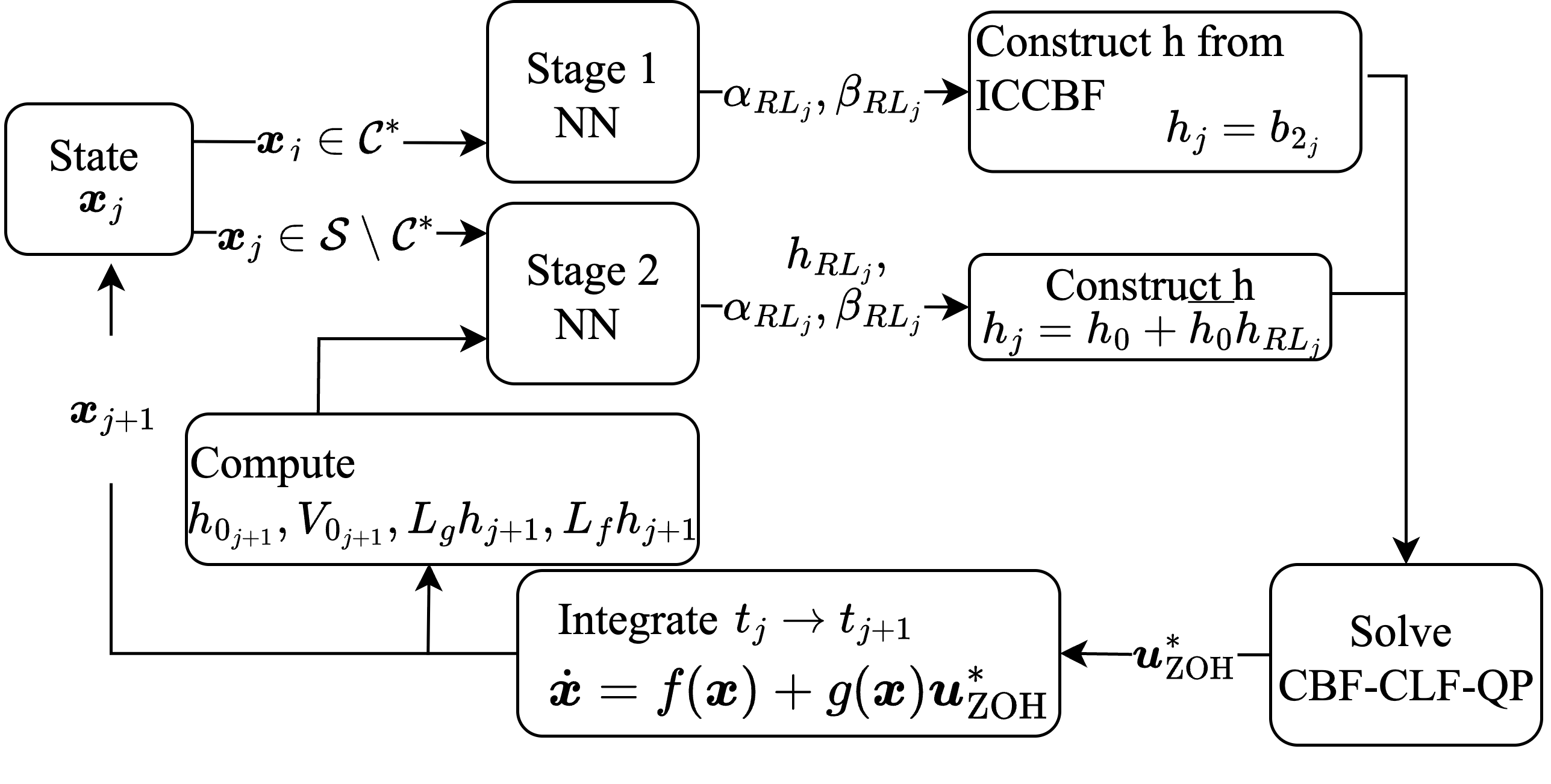}
    \caption{Flow chart of the two-stage RL framework.}
    \label{fig:overallmethod}
  \end{subfigure}\hfill
  \begin{subfigure}[t]{0.5\textwidth}
    \centering
    \includegraphics[width=1.1\linewidth]{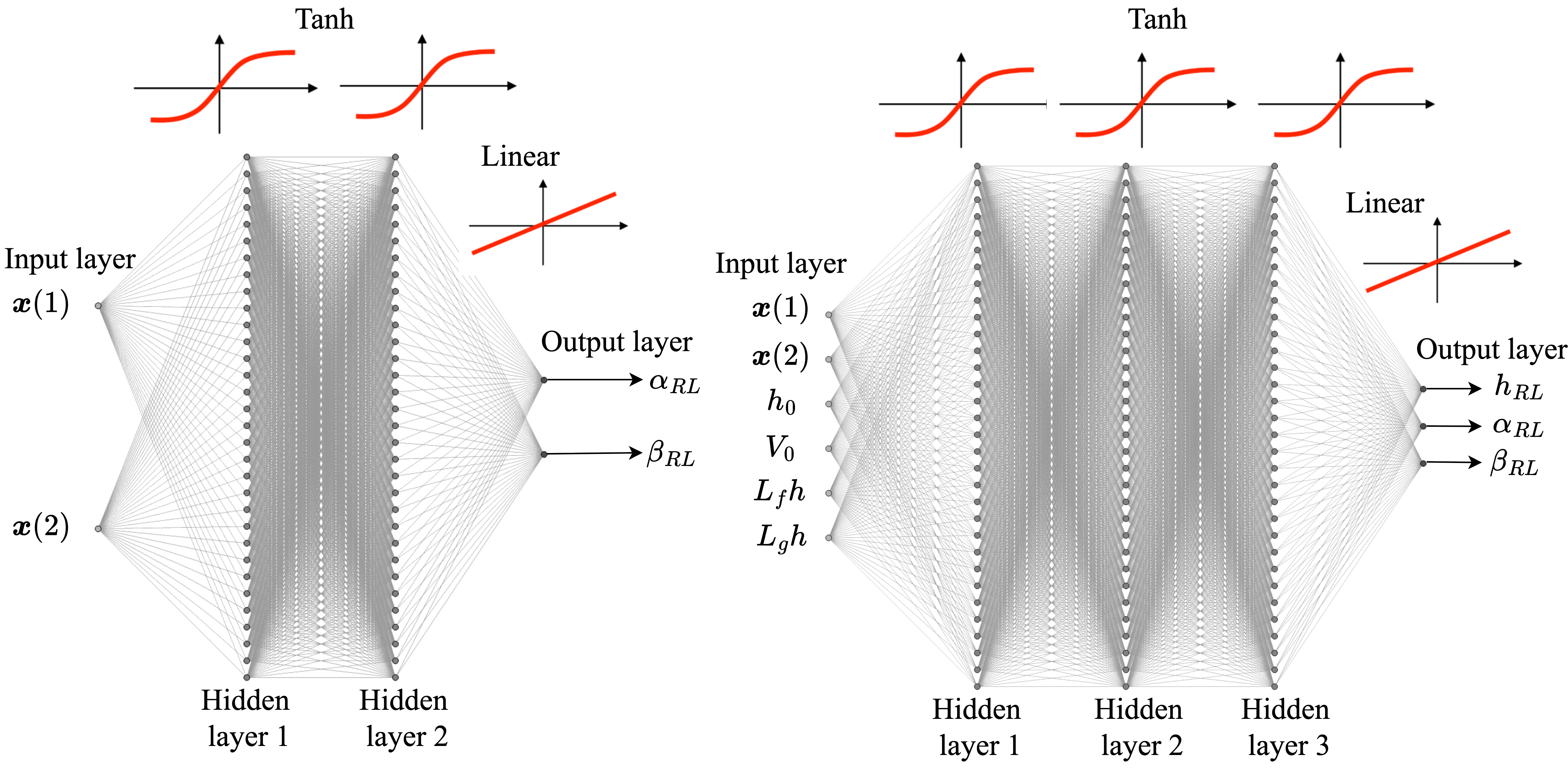}
    \caption{Actor-network architectures for Stage~1 (left) and Stage~2 (right), modelled for the Cruise Control problem (Sec.~\ref{ccproblemintro}).}
    \label{fig:networks}
  \end{subfigure}
  \caption{Overview of the framework and corresponding actor-network designs.}
  \label{fig:overview-and-nets}
\end{figure*}

As mentioned in the Introduction, Stage~1 mitigates greediness within \(C^*\) by learning state-dependent class-$\mathcal{K}_\infty$ functions for ICCBF, improving long-horizon efficiency while preserving invariance. 
Stage~2 enlarges the effective operating region by learning a residual barrier on \(\mathcal{R}\) that certifies recoverability for states otherwise discarded by the baseline ICCBF. An overall schematic of how these two learnt controllers work together is given in Figure \ref{fig:overallmethod}, and example network architectures for both stages are shown in Figure \ref{fig:networks}. Note that Stage 1 is used for the current state is in the inner safe set, i.e, $x \in  \mathcal{C}^*$. Stage 2 is used when $x$ is still in the safe region but not within the inner safe set, i.e, $x \in  \mathcal{S} \setminus \mathcal{C}^*$.
Together, the two stages yield CBFs that are less greedy within \(C^*\) and work beyond \(C^*\).



Consider a control-affine system given in Eq. \eqref{xdoteq} with admissible inputs $u \in \mathcal{U}$. The methodology consists of the following common steps for both RL stages. 
\begin{enumerate}
    \item At each decision step, an RL policy outputs CBF-related parameters. The parameters themselves are different for Stage 1 and Stage 2, but are used to generate a CBF-QP in each case.
    \item The generated CBF-(CLF)-QP is solved to compute the control input $u^*$ to retain the invariance of the CBF. 
    \item The control is applied under zero-order-hold (ZOH) for a specified time interval while the system dynamics are forward-propagated.
    \item  RL agent is given a reward based on fuel consumption (to minimize), violations of the original safety constraints after propagation, and optionally other task-specific performance metrics. 
\end{enumerate}

\subsection{Stage 1: RL-tuned ICCBF}\label{subsec:scaled-ICCBF}

This stage preserves the analytical structure of an ICCBF and performs adaptive tuning of class-$\mathcal{K}_\infty$ parameters $\alpha_{RL}$ and $\beta_{RL}$ to reduce greediness of the trajectories. Firstly, let $b_1, b_2$ denote the first two ICCBF layers constructed from a nominal safety function $h_0$. Then, the barrier constraint is defined as  
\begin{equation}\label{hofmethod1}
    h(\boldsymbol{x}) = b_2(\boldsymbol{x}).
\end{equation}
During training the initial states $\boldsymbol{x}_0$ are drawn from the inner safe set $\mathcal{C}^*$ defined in Eq. \eqref{cstar}. A PPO policy learns $\alpha_{RL}(\boldsymbol{x}) \in [\alpha_{\min}, \alpha_{\max}]$ and $\beta_{RL}(\boldsymbol{x}) \in [\beta_{\min}, \beta_{\max}]$ that replace the fixed parameters in the standard ICCBF formulation.
The following QP is then solved to determine the control needed to maintain invariance:
\begin{equation}
\begin{aligned}
& \argmin_{u \in \mathbb{R}^2, \; \delta, \gamma}
& & \frac{1}{2} u^T u + p_1 \delta + p_2 \gamma \\
&\text{subject to}
&& L_f h(\boldsymbol{x}) + L_g h(\boldsymbol{x}) \, u \geq -(\alpha_{RL}(\boldsymbol{x}) + \gamma) h(\boldsymbol{x}), \\
& & & L_f V(\boldsymbol{x}) + L_g V(\boldsymbol{x}) \, u \leq -\beta_{RL}(\boldsymbol{x}) V(\boldsymbol{x}) + \delta, \\
&&& \|u\| \leq u_{\max}
\end{aligned}
\label{QPm1}
\end{equation}
where $p_1$ and $p_2$ are user-set penalty coefficients for the slack variables $\delta$ and $\gamma$, respectively.

\subsection{Stage 2: RL–residual CBF}\label{subsec:learned-CBF}

It can be noted that the ICCBF construction utilizes worst-case analysis, as the definitions of $b_0,..., b_N$ in Eq.~\eqref{bees} utilize the infimum over $\mathcal{U}$. Furthermore, the use of conservative class-$\mathcal{K}_{\infty}$ functions may also exclude states from $C^*$ that could potentially be rendered safe under trajectory-aware control policies. However, trajectory-aware control policies operating over finite horizons can potentially employ strategies unavailable to ICCBFs, such as accepting controlled constraint relaxation at time  $t$ if it provides better safety margins at $t+ \Delta t$ and learning state-dependent $\mathcal{K}_{\infty}$ functions rather than 
fixed ones.  The fundamental gap between instantaneous worst-case analysis (ICCBF) and finite-horizon trajectory optimization creates the theoretical possibility for recovery of some states in $\mathcal{R} = \mathcal{S} \setminus C^*$. 

While not all states in $\mathcal{R}$ are necessarily recoverable due to genuine control limitations, the conservative nature of the ICCBF construction suggests that some excluded states may admit safe trajectories under non-greedy policies. This provides the theoretical motivation for Stage 2 of the proposed framework. 
Thus, in this section, RL attempts to identify and recover such abandoned-but-potentially-safe states through the learned residual barrier function $h_{RL}(x)$.

Here, PPO learns $h_{RL}$, $\alpha$, and $\beta$, where the barrier function is defined as:
\begin{equation}\label{hofmethod2}
    h(\boldsymbol{x})= h_{0}(\boldsymbol{x}) + \overline{h}_0 h_{RL}(\boldsymbol{x})
\end{equation}
where $\overline{h}_0$ is the mean value of $h_{0}$ over the initial distribution discussed in Section \ref{results_benchmark}, providing appropriate scaling for the learned residual term. During training and evaluation, the initial states $\boldsymbol{x}_0$ are drawn from the residual set $\mathcal{R}$ in this stage.  During training and evaluation, initial states are sampled uniformly from $\mathcal{R}$ as the objective of Stage 2 is the recovery of some of the states abandoned by ICCBF due to conservatism. The quadratic program formulation remains identical to Eq. \eqref{QPm1}, but now employs the learned barrier function $h(\boldsymbol{x})$ from Eq. \eqref{hofmethod2}.

The gradient $ \frac{\partial h}{\partial \boldsymbol{x}}$ is required to derive the lie derivatives in the QP formulation, and is given by 
\begin{equation}
    \frac{\partial h}{\partial \boldsymbol{x}} = \frac{\partial h_0}{\partial \boldsymbol{x}} + \overline{h}_0\frac{\partial h_{RL}}{\partial \boldsymbol{x}}
\end{equation}
where $\frac{\partial h_{RL}}{\partial \boldsymbol{x}}$ is obtained via automatic differentiation \cite{paszke2017automatic}.

\subsection{Actor Input}
For Stage 1, the actor input state $\boldsymbol{S}_j$ is set to be $\boldsymbol{x}_j$.   For Stage 2, the observation space is extended to $\boldsymbol{S}_j = [\boldsymbol{x}_j, L_g h(\boldsymbol{x}_j), L_f h(\boldsymbol{x}_j), h_0(\boldsymbol{x}_j), V(\boldsymbol{x}_j)]$, as this combination was noted to yield faster convergence.  In both cases, the actor inputs are normalised to lie in $[-1,1]$, promoting numerical stability and consistent feature scaling~\cite{WIJAYATUNGA2025109996}. Thus
\begin{equation}\label{actorin}
\boldsymbol{S}_j^*=2\left(\frac{\boldsymbol{S}_j-\boldsymbol{S}_{\min}}{\boldsymbol{S}_{m a x}-\boldsymbol{S}_{\min}}\right)-1.
\end{equation} 
where $\boldsymbol{S}_{\min}$ and $\boldsymbol{S}_{\max}$  are the vectors of minimum and maximum values of the components in $\boldsymbol{S}$, respectively.

\subsection{Reward}
 For both stages, the RL reward at step $j$ ($R_j$) penalizes safety violations and control effort:
\begin{equation}\label{method1Penalty}
R_j = -c_h \max(0, -h_0(\boldsymbol{x}_j)) - c_u \|\boldsymbol{u}_j\|
\end{equation}
where $c_h$ and $c_u$ are user-set cost scaling factors. This reward structure encourages the RL agent to a policy that maintains safety while minimizing control effort, effectively learning when to be conservative versus aggressive based on the current system state and mission phase.

Note that for both cases, the RL actor's inputs, outputs, and rewards are summarized in Table \ref{t1}, and the state propagation and reward generation in the RL environment are discussed in Algorithm \ref{GE}.
\begin{table}[hbt!]
\centering 
\renewcommand{\arraystretch}{1.3}
\begin{tabular}{lll}\hline 
Parameter & Step 1 & Step 2 \\ \hline
Initial State Set & $\boldsymbol{x}\in \mathcal{C}^*$&  $\boldsymbol{x}\in \mathcal{S} \setminus \mathcal{C}^*$  \\
Actor Input State $\boldsymbol{S}_j$   & $ \boldsymbol{x}_j $ &  $[\boldsymbol{x}_j, L_g h(\boldsymbol{x}_j) , L_f h(\boldsymbol{x}_j) ,$ \\
& & $ h_0(\boldsymbol{x}_j), V(\boldsymbol{x}_j)]$\\
Actor Output/Action&  $\alpha_{RL},\beta_{RL}$ & $h_{RL}, \alpha_{RL},\beta_{RL}$ \\
Reward              &   \multicolumn{2}{l}{Given in Eq.~\eqref{method1Penalty}}\\
State Transition    &   \multicolumn{2}{l}{Given in Algorithm \ref{GE}.}\\ \hline 
\end{tabular} 
\caption{{An overview of the  RL setup}}
\label{t1}
\end{table}

\begin{algorithm*}[hbt!]
\renewcommand{\arraystretch}{1.3}
\caption{RL environment}
\label{GE}
\textbf{Input:} PPO actions $\mathbf{a}_j=[\alpha_{RL,j},\beta_{RL,j}]$ (Stage 1) or $[h_{\mathrm{RL},j},\alpha_{RL,j},\beta_{RL,j}]$ (Stage 2); state $\boldsymbol{x}_j$; safety function $h(\boldsymbol{x})$; current time $t$; time step $\Delta t$; final time $t_f$; ICCBF $b_2(\boldsymbol{x})$.

\begin{algorithmic}[1]
\State Compute $h(\boldsymbol{x}_j)$ using Eq.~\eqref{hofmethod1} (for Stage 1) or Eq.~\eqref{hofmethod2} (for Stage 2)
\State Compute $\partial h/\partial \boldsymbol{x}$ and the corresponding Lie derivatives.
\State Compute the ZOH control $\boldsymbol{u}_j$ by solving the QCQP in~\eqref{QPm1}.
\State Forward propagate the system on $[t,\,t+\Delta t]$ under the dynamics and $\boldsymbol{u}_j$ to obtain $\boldsymbol{x}_{j+1}$.
\State Update $t \gets t+\Delta t$.
\State Compute the reward $R_{j+1}$ using~\eqref{method1Penalty}.
\State Calculate the propagated actor input state:
\[
\boldsymbol{S}_{j+1} =
\begin{cases}
[\boldsymbol{x}_{j+1}], & \text{if Stage 1},\\
[\boldsymbol{x}_{j+1},\, L_g h(\boldsymbol{x}_{j+1}),\, L_f h(\boldsymbol{x}_{j+1}),\, h_0(\boldsymbol{x}_{j+1}),\, V(\boldsymbol{x}_{j+1})], & \text{if Stage 2}.
\end{cases}
\]
\State Scale $\boldsymbol{S}_{j+1}$ to obtain $\boldsymbol{S}_{j+1}^\ast$ using~\eqref{actorin}.
\State \textbf{Output:} $\boldsymbol{S}_{j+1}^\ast,\ R_{j+1}$.
\If{$t = t_f$}
  \State \textbf{Stop.}
\Else
  \State \textbf{Update:} $j \gets j+1$ and repeat from Step~1.
\EndIf
\end{algorithmic}
\end{algorithm*}

\section{Results: Benchmark Problems}\label{results_benchmark}

This section outlines the problem dynamics, implementation details, and results for two benchmark test cases: cruise control and spacecraft docking with a rotating target, adapted from \cite{agrawal2021safe}.

\subsection{Test Cases}
\subsubsection{Cruise Control Problem}\label{ccproblemintro}

\begin{figure}[hbt!]
    \centering\includegraphics[width=\linewidth]{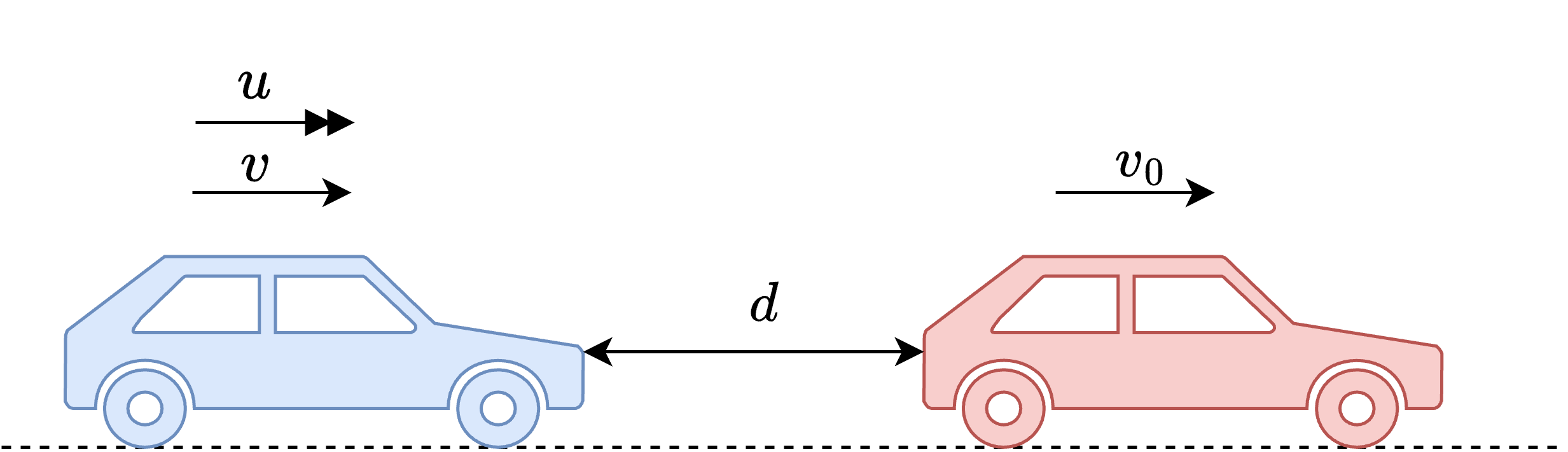}
    \caption{The Cruise Control Problem}
    \label{fig:cc}
\end{figure}

Shown in Figure \ref{fig:cc} and discussed in \cite{7040372} and \cite{agrawal2021safe}, this problem considers a point-mass model of a vehicle moving along a straight line. A following vehicle trails a lead vehicle at a distance $d$, with the lead vehicle travelling at a known constant speed $v_0$. The objective is to design a controller that drives the following vehicle to the speed limit $v_{\text{max}} = \SI{24}{\meter\per\second}$ while ensuring collision avoidance. The collision avoidance safety constraint is specified as $d \geq 1.8v$, and the CLF constraint that drives the following vehicle to the speed limit
is defined as $V(x)  =(v - v_{\text{max}})^2$. {No drag is considered in this preliminary cruise control scenario.} Defining the state vector as $x = [d,\, v]^T$, the dynamical model is
\begin{equation}
\begin{bmatrix}
\dot{d} \\
\dot{v}
\end{bmatrix}
=
\begin{bmatrix}
v_0 - v \\
-\tfrac{F(v)}{m}
\end{bmatrix}
+
\begin{bmatrix}
0 \\
g_0
\end{bmatrix} u,
\qquad
\mathcal{U} = \{u : |u| \leq \SI{0.25}{\newton}\},
\end{equation}
where $u$ is the control input. The resistive force $F(v)$ is modelled as $F(v) = f_0 + f_1 v + f_2 v^2$, $m$ is the vehicle mass, and $g_0$ is the gravitational acceleration.  The safe set $\mathcal{S}$ is then defined as
\begin{equation}
\mathcal{S} = \{x \in \mathcal{X} \;|\; h_0(x) = x_1 - 1.8x_2 \geq 0\}.
\end{equation}
\subsubsection{Spacecraft Rendezvous with Rotating Target}

Discussed in \cite{agrawal2021safe} and \cite{5991151}, this test case is an autonomous rendezvous scenario between a chaser spacecraft and a target body, both modelled as point masses. The target is represented as a point on a disk of radius $\rho = 2.4 \,\text{m}$, rotating with a constant angular velocity $\omega = 0.6^\circ/\text{s}$ relative to the Local-Vertical Local-Horizontal (LVLH) frame. The objective is to determine the thrust needed to bring the chaser spacecraft from an initial range of 100~m to within 3~m of the target.  

\begin{figure}[hbt!]
\centering\includegraphics[width=\linewidth]{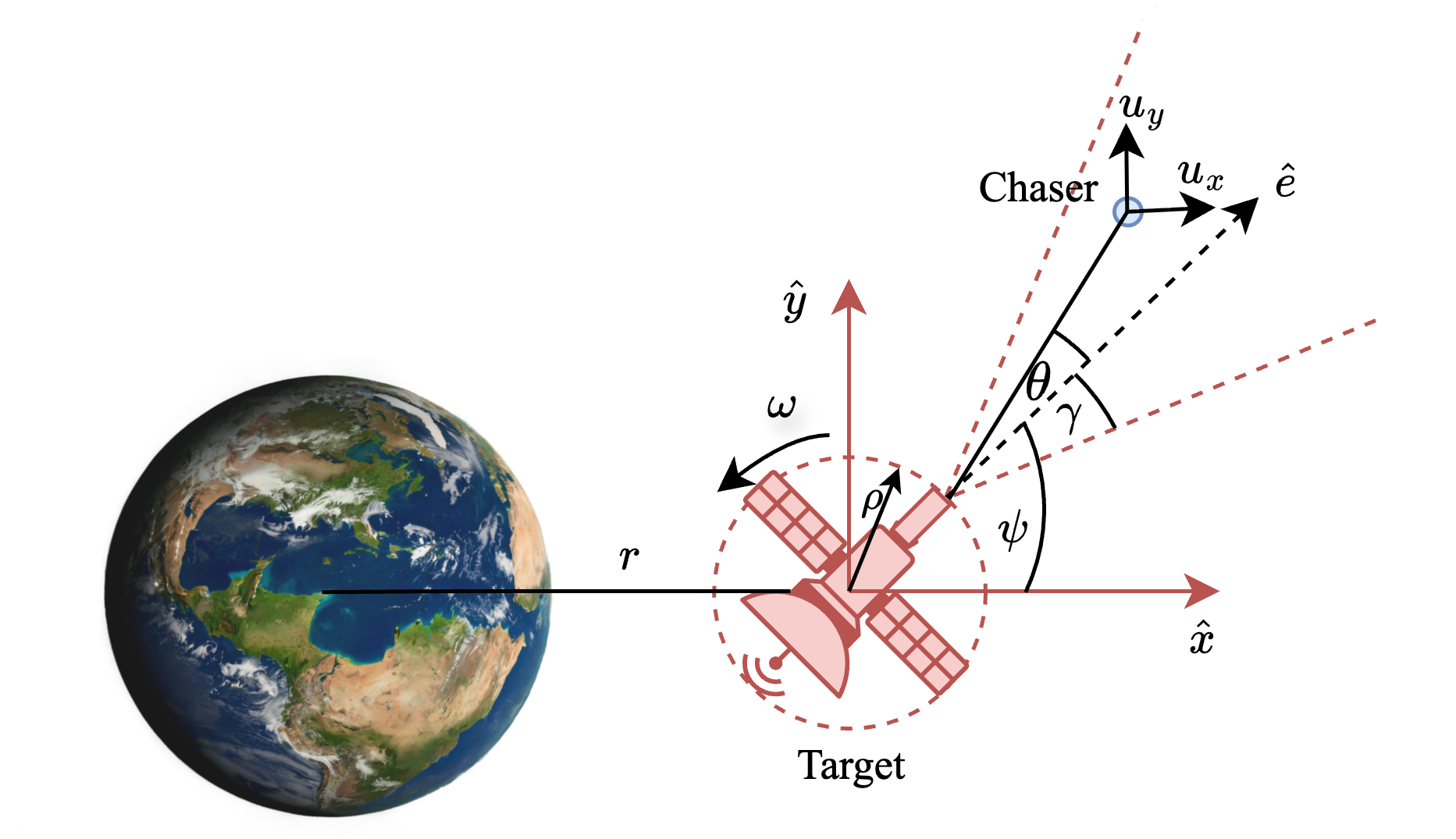}
    \caption{Spacecraft Docking Problem}
    \label{dockingproblemFig}
\end{figure}

A line-of-sight (LOS) safety constraint is imposed in this test case, requiring that the chaser’s relative position remain inside a cone of half-angle $\gamma = 10^\circ$ aligned with the docking axis at all times. The system state is defined as $x \in \mathbb{R}^5$, comprising the relative position $(p_x, p_y)$, relative velocity $(v_x, v_y)$, and the docking port angle $\psi$ as shown in Figure \ref{dockingproblemFig}.  

The system dynamics are then given as 
\begin{equation}
\frac{d}{dt}
\begin{bmatrix}
p_x \\ p_y \\ v_x \\ v_y \\ \psi
\end{bmatrix}
=
\begin{bmatrix}
v_x \\
v_y \\
n^2 p_x + 2n v_y + \tfrac{\mu}{r^2} - \tfrac{\mu(r+p_x)}{r_c^3} \\
n^2 p_y - 2n v_x - \tfrac{\mu p_y}{r_c^3} \\
\omega
\end{bmatrix}
+
\frac{1}{m_c}
\begin{bmatrix}
0 \\ 0 \\ u_x \\ u_y \\ 0
\end{bmatrix},
\label{eq:dyn}
\end{equation}
where $r_c = \sqrt{p_x^2 + p_y^2}$ is the relative distance between the vehicles, $r = 6771 \,\text{km}$ is the orbital radius of the target, and $\mu = 398{,}600 \,\text{km}^3/\text{s}^2$ is the Earth’s gravitational parameter. The mean motion of the target is $n = \sqrt{\mu/r^3}$  and $m_c = 1000 \,\text{kg}$ is the chaser spacecraft mass.   The control inputs $\boldsymbol{u} = [u_x, u_y]^T$ represent the propulsive forces and are bounded such that $\parallel \boldsymbol{u}\parallel  \leq 0.25 \,\text{kN}$.   The LOS constraint is expressed as $h_0(x) \geq 0$, where  
\begin{align}
h_0(x) &= \cos\theta - \cos\gamma \\
     &= \frac{\vec{r}_{c-p} \cdot \hat{e}}{\|\vec{r}_{c-p}\|} - \cos(\gamma).
\end{align}
Here, $\vec{r}_{c-p} = \big(p_x - \rho \cos\psi, \; p_y - \rho \sin\psi\big)^T$ is the position vector of the chaser relative to the docking port, and $\hat{e} = (\cos\psi, \sin\psi)^T$ is the docking axis unit vector.  Note that this $h_0(x)$ has relative degree challenges as $L_g h_0(x) = 0$. 

To guide the chaser to the docking port the following CLF is used 
\begin{equation}
V(x) = 
\left( v_x + \frac{p_x - \rho \cos\psi}{10} \right)^2
+ \left( v_y + \frac{p_y - \rho \sin\psi}{10} \right)^2.
\end{equation}

\subsection{Implementation Details}

\subsubsection{Training Configuration}
When solving the benchmark cases, the PPO algorithm was implemented using Stable Baselines \cite{raffin2021stable}. Both stages utilize policy networks with 64 neurons each, employing hyperbolic tangent (Tanh) activation functions throughout. The complete hyperparameter configuration for the two baseline test cases for Stages 1 and 2 is provided in Table~\ref{tab:hyperparameters}. These parameters were chosen empirically by running small training batches with varying settings in a simple parameter sweep

Stage 2 training requires computing $\partial h/\partial \boldsymbol{x}$ with access to the policy during training. Due to the Stable Baselines architecture, this prevents multi-threading, necessitating single-thread environments for all Stage 2 simulations.

\begin{table*}[hbt!]
\centering
\caption{PPO Hyperparameters for Benchmark Problems}
\label{tab:hyperparameters}
\renewcommand{\arraystretch}{1.0}
\begin{tabular}{|l|c|c|c|c|}
\hline
{Parameter} & \multicolumn{2}{c|}{{Cruise Control}} & \multicolumn{2}{c|}{{Docking}} \\
\cline{2-5}
 & {Stage 1} & {Stage 2} & {Stage 1} & {Stage 2} \\ \hline 
Learning Rate & $10^{-3}$ (const.) & $10^{-4}$ (decay.) & $10^{-3}$ (const.) & $10^{-4}$ (decay.) \\
Batch Size & 64 & 256 & 64 & 256 \\
Rollout Steps & 1,280 & 1,280 & 2,560 & 2,560 \\
Training Epochs & 10 & 10 & 10 & 10 \\
Discount Factor ($\gamma$) & 0.95 & 0.999 & 0.95 & 0.999 \\
GAE Lambda ($\lambda$) & 0.99 & 0.99 & 0.99 & 0.99 \\
Clip Range & 0.2 & 0.2 & 0.2 & 0.2 \\
Entropy Coefficient & 0.01 & 0.01 & 0.01 & 0.01 \\
Initial Std. Dev. & 0.2 & 0.2 & 0.2 & 0.2 \\
SDE & Yes & No & Yes & No \\
Network Layers & 4 & 4 & 4 & 4 \\
Neurons per Layer & 64 & 64 & 64 & 64 \\
Activation Function & Tanh & Tanh & Tanh & Tanh \\
Parallel Environments & 8 & 1 & 8 & 1 \\
Total Episodes & $10^6$ & $10^4$ & $10^6$ & $10^4$ \\
\hline 
\end{tabular}
\end{table*}
For both stages, the runtime computational overhead is minimal: Stage 1 requires only forward neural network evaluation ($\approx 0.1$ ms), in addition to standard QP solving, while Stage 2 adds automatic differentiation ($\approx 0.3$ ms additional). Training requires 12-48 hours on standard hardware but is performed offline.

\subsubsection{Episode Termination and Evaluation}
Episodes terminate upon reaching the maximum simulation time ($t_f = 20$s for cruise control, $t_f = 50$s for docking, with time steps $\Delta t = 0.1$s and $0.5$s respectively). For the docking scenario, episodes may terminate early when successful docking is achieved ($V < 5 \times 10^{-5}$), resulting in approximate episode lengths of 200 steps (cruise control) and 100 steps (docking).

Model evaluation is performed after each rollout buffer completion using 10 deterministic episodes on a separate evaluation environment.  The initial conditions for the evaluation episodes are uniformly sampled from the feasible region $\mathcal{C}^*$ for Stage 1 and from the broader region $\mathcal{R}$ for Stage 2.
The best-performing model, based on cumulative reward across evaluation episodes, is saved for final testing. Training utilizes the Adam optimizer with advantage normalization enabled for improved convergence stability.

\subsubsection{Initial State Distributions}
\quad

\textbf{Cruise Control:}
The initial state space is discretized as 
$\mathcal{G} = \{(x_1^i, x_2^j)\}$ where $x_1^i \in \{0, 10, 20, \ldots, 120\}$ and  $x_2^j \in \{0, 1, 2, \ldots, v_{\max}\}$.

For ICCBF and Stage 1 methods, initial states are sampled from the inner safe set:
\begin{equation}\label{Ddist}
\mathcal{D} = \{\boldsymbol{x} \in \mathcal{G} : h(\boldsymbol{x}) \geq 0 \land b_1(\boldsymbol{x}) \geq 0 \land b_2(\boldsymbol{x}) \geq 0\} \subset \mathcal{C}^*
\end{equation}
For Stage 2 evaluation, states are sampled from the expanded feasible region:
\begin{equation}\label{Edist}
\mathcal{E} = \{\boldsymbol{x} \in \mathcal{G} : h(\boldsymbol{x}) \geq 0 \land (b_1(\boldsymbol{x}) < 0 \lor b_2(\boldsymbol{x}) < 0)\} \subset \mathcal{R}
\end{equation}

\textbf{Docking:}
Initial conditions are geometrically distributed within the line-of-sight cone. States are defined as $\boldsymbol{x} = [x_1, x_2, v_x, v_y, \psi]^T$ where $x_1 = 500$~m represents the starting standoff distance, $x_2 = x_1/\tan(\theta_i)$ determines the lateral position, velocities are initialized to zero ($v_x = v_y = 0$), and the docking port angle is $\psi = 0$. The angular positions $\theta_i$ are uniformly distributed within the cone boundary:
$$\theta_i \in \{\theta_j : \theta_j = -\gamma + \frac{2\gamma j}{99}, j = 0, 1, \ldots, 99\}$$
ensuring 100 equally-spaced initial conditions along the approach cone edge.

For ICCBF and Stage 1 methods, the feasible initial states are given by $\mathcal{D} \subset \mathcal{C}^*$ as defined in Eq. \eqref{Ddist}, while Stage 2 operates on initial states from $\mathcal{E} \subset \mathcal{S} \setminus \mathcal{C}^*$ as given in Eq. \eqref{Edist}. 

\subsection{Numerical Results}
\subsubsection{Cruise Control}
\begin{figure*}[hbt!]
    \centering
\includegraphics[width=\linewidth]{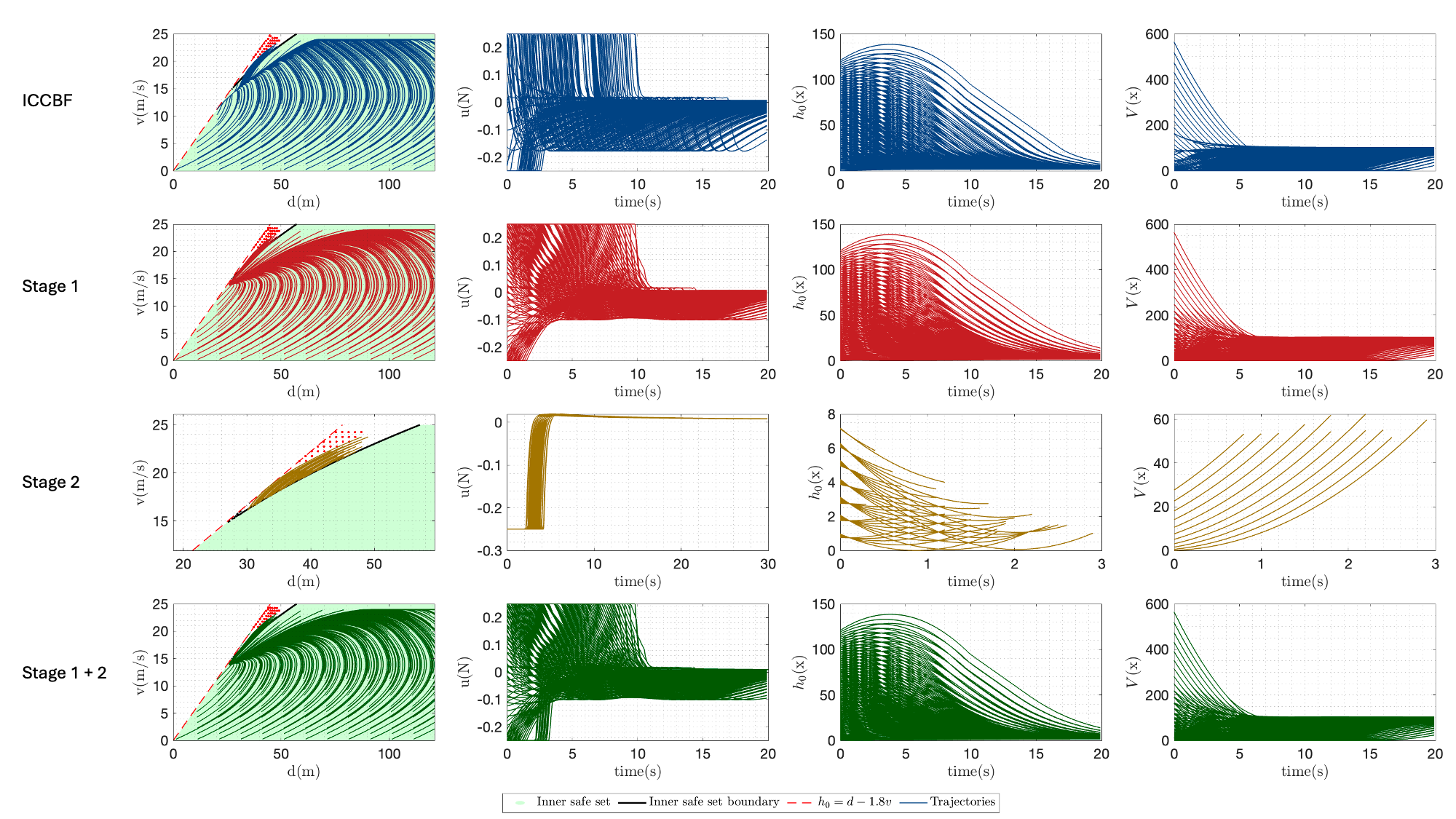}
    \caption{Cruise Control Results: Comparative analysis of trajectory performance, control profiles, safety constraints ($h_0$), and CLF evolution. Failed cases are shown in red.}
    \label{fig:ccres}
\end{figure*}

\begin{figure}[hbt!]
    \centering
\includegraphics[width=\linewidth]{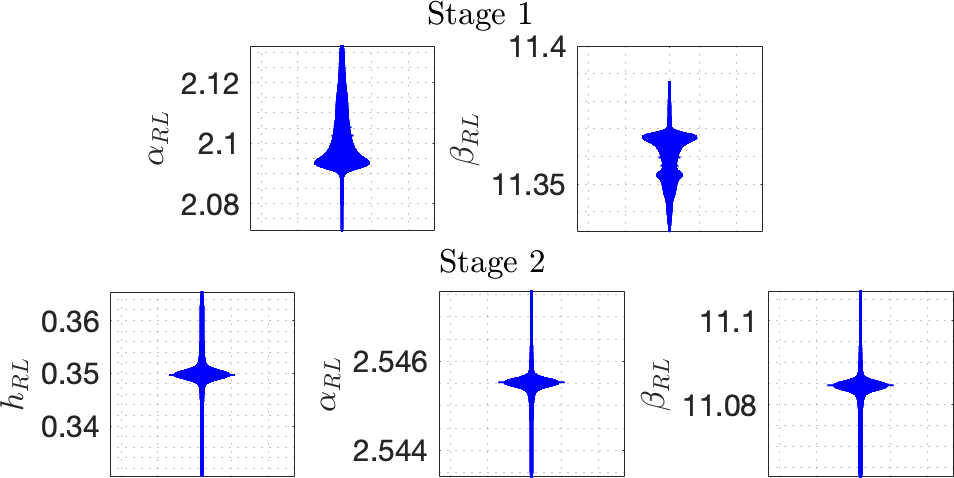}
    \caption{Neural network outputs for the cruise control case}
    \label{fig:ccNNout}
\end{figure}

\begin{figure}[hbt!]
    \centering
\includegraphics[width=\linewidth]{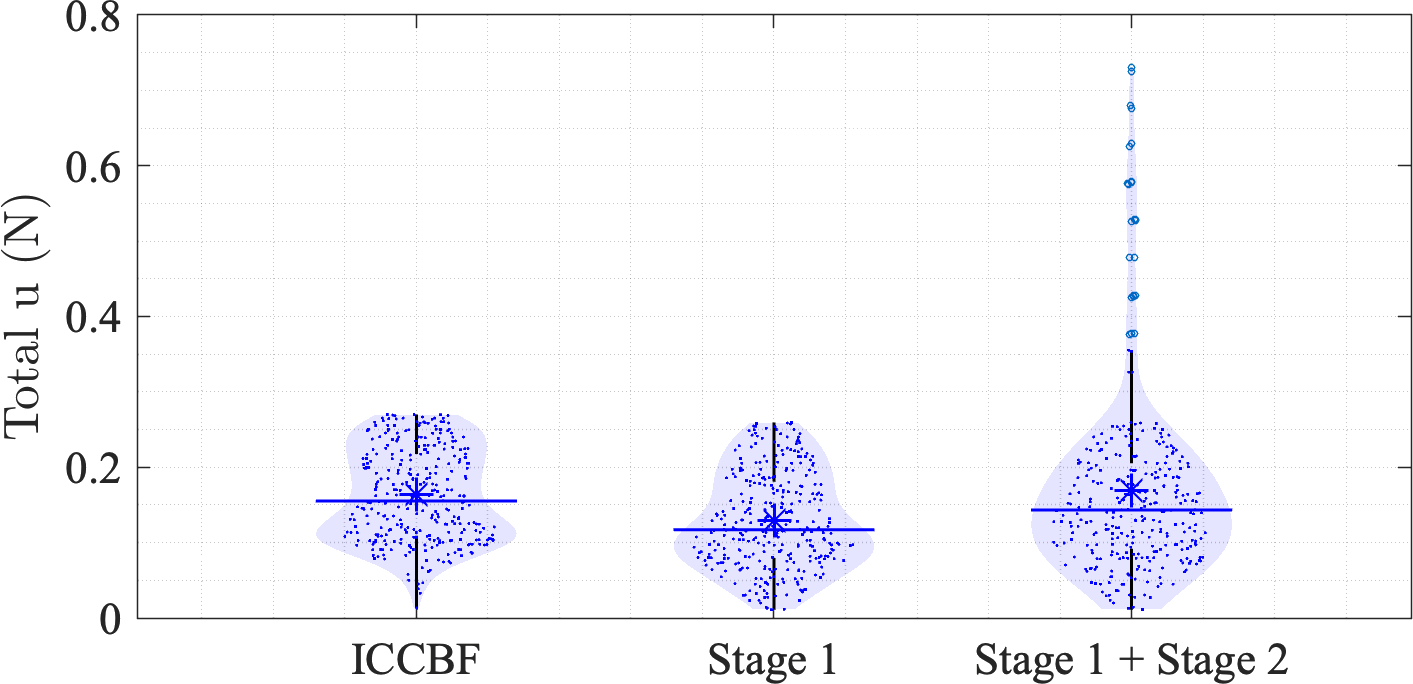}
    \caption{Control effort distribution for cruise control scenarios across different strategies.}
    \label{fig:CCviolin}
\end{figure}

The cruise control problem demonstrates the effectiveness of the two-stage approach through comparison of three control strategies: baseline ICCBF, Stage 1 (RL-tuned ICCBF), and Stage 2 (RL residual CBF) framework. Figure \ref{fig:ccres} presents the comparative trajectory analysis, control profiles, safety function evolution ($h_0 = d-1.8v$), and CLF behavior across the three approaches and the combination of Stage 1 and 2. Initial states where the safety cannot be satisfied under the control barrier strategy utilized are shown as red dots. The total control magnitude distributions for each case are visualized in Figure \ref{fig:CCviolin}, with their statistical summaries provided in Table \ref{tableCC}. The distributions of the action states for both Stages are given in Figure \ref{fig:ccNNout}.

\begin{table*}[hbt!]
\centering
\caption{Performance Comparison of CBF Strategies Across the Benchmarking Scenarios. (The Median Change is with respect to the ICCBF median fuel consumption)}
\label{tableCC}
\begin{tabular}{lccccc}
\hline
Strategy & Successes & Success (\%) & \multicolumn{3}{c}{Fuel Consumption [N]} \\
\cline{4-6}
            &                         &             & \textbf{$\mu \pm \sigma$} & [Q$_1$, Q$_2$, Q$_3$, P$_{99}$] & $\Delta Q_2$ (\%) \\
\hline

\multicolumn{6}{l}{\textbf{Cruise Control}} \\
\hline
ICCBF       & 279/325 & 85.8\% & 0.16 $\pm$ 0.06 & [0.11, 0.16, 0.22, 0.27] & 0\% \\
Stage 1     & 279/325 & 85.8\% & 0.13 $\pm$ 0.06 & [0.08, 0.12, 0.18, 0.26] & -25.0\% \\
Stage 1+2   & 299/325 & \textbf{92.0\%} & 0.17 $\pm$ 0.13 & [0.09, 0.14, 0.21, 0.68] & \textbf{-12.5\%} \\
\hline

\multicolumn{6}{l}{\textbf{Docking}} \\
\hline
ICCBF       & 88/100  & 88.0\% & 0.74 $\pm$ 0.04 & [0.71, 0.72, 0.77, 0.82] & 0\% \\
Stage 1     & 88/100  & 88.0\% & 0.66 $\pm$ 0.04 & [0.64, 0.64, 0.66, 0.78] & -11.1\% \\
Stage 1+2   & 95/100  & \textbf{95.0\%} & 0.69 $\pm$ 0.11 & [0.64, 0.64, 0.66, 1.00] & \textbf{-11.1\%} \\
\hline
\end{tabular}
\end{table*}

The results demonstrate significant improvements across multiple performance metrics. Stage 1 achieves superior fuel efficiency compared to baseline ICCBF, with a 25\% reduction in median total control magnitude while maintaining safety guarantees within $\mathcal{C}^*$. The two-stage approach significantly expands the feasible operating region, as both ICCBF and Stage 1 are limited to $\mathcal{C}^* \subset \mathcal{S}$, leaving states in $\mathcal{S} \setminus \mathcal{C}^*$ unaddressed, while the combined Stage 1+2 framework can handle more of the broader safe set $\mathcal{S}$. Within the 325 evaluation points, Stage 1 and ICCBF exhibit 46 failed cases (85.8\% success rate) when starting state is in $\mathcal{S} \setminus \mathcal{C}^*$, but integration of Stage 1+2 reduces this to 26 failures (92.0\% success rate), representing a significant improvement in success rate.  Note that there are no failure cases when  $x_0 \in \mathcal{C}^*$, as expected. When $ x_0 \in \mathcal{S} \setminus \mathcal{C}^*$, ICCBF theory dictates that thrust is insufficient to remain in the safe region $\mathcal{S}$, however, the combination of Stage 1 and Stage 2 can recover some of these points and retain the trajectories in the safe region.  The combined Stage 1 and 2 has a higher median fuel consumption than Stage 1, representing the additional control required to recover trajectories from $\mathcal{S} \setminus \mathcal{C}^*$.  Note that, once training is complete, each timestep of the control generation is extremely fast and can be executed in real time, since it only requires solving a single QP. The resulting average computation times for ICCBF, Stage~1, and Stage~2 are \(0.0036 \pm 0.0025\)~ms, \(0.0034 \pm 0.0025\)~ms, and \(0.0036 \pm 0.0027\)~ms, respectively.

\subsubsection{Rendezvous Results}
\begin{figure*}[hbt!]
    \centering
\includegraphics[width=\linewidth]{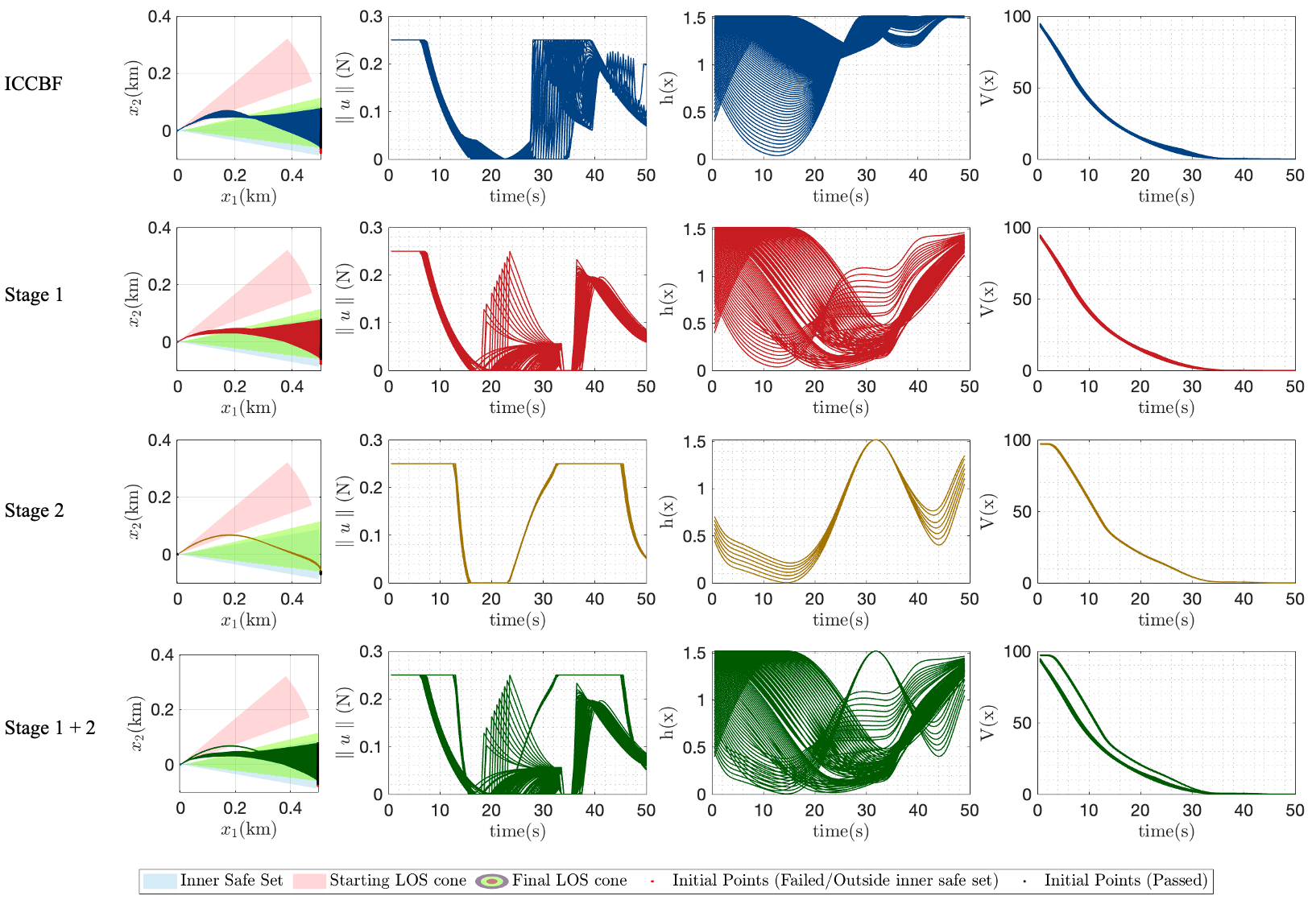}
    \caption{Rendezvous and Docking Results: Comparative analysis of trajectory performance, control profiles, safety constraints ($h_0$), and CLF evolution across different control strategies. Failed cases are shown in red.}
    \label{fig:dockingres}
\end{figure*}

\begin{figure}[hbt!]
    \centering
    \includegraphics[width=\linewidth]{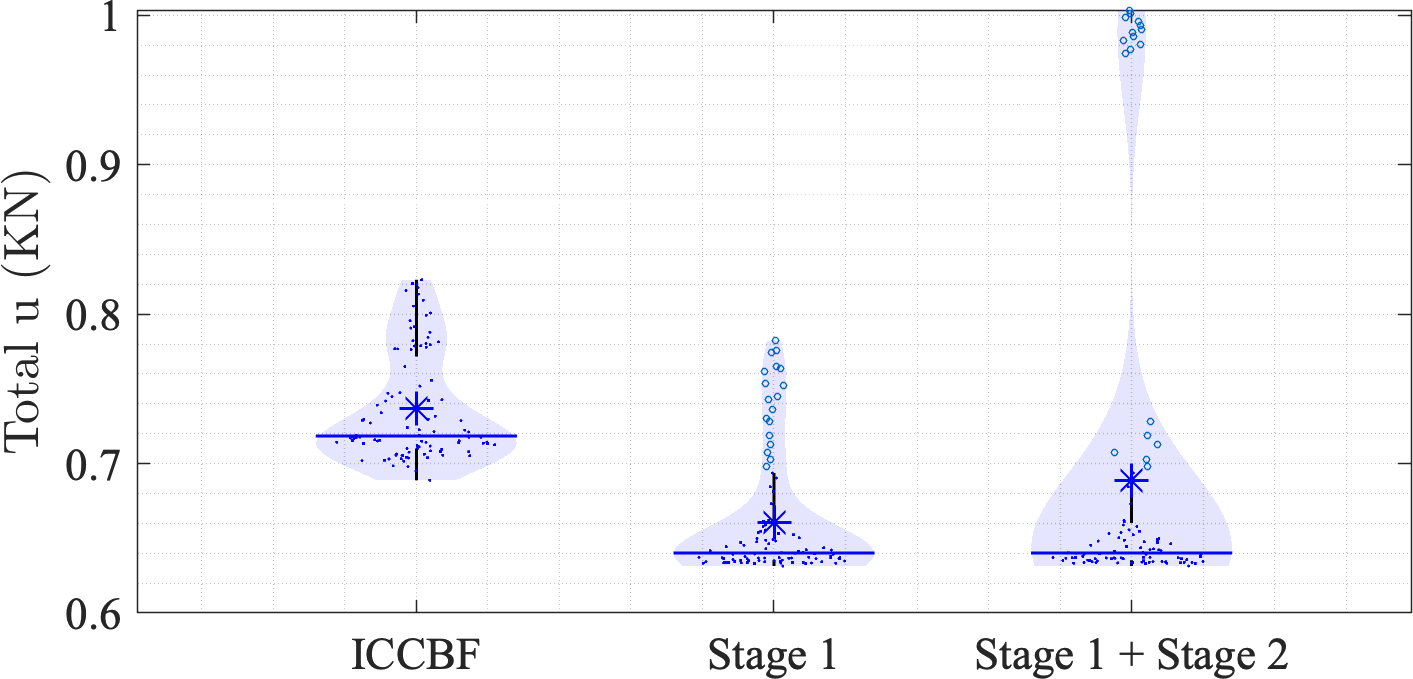}
    \caption{Control effort distribution for rendezvous and docking scenarios across different strategies.}
    \label{fig:dockingcontrolthrust}
\end{figure}

This test case evaluates the performance of RL control strategies for spacecraft docking with a rotating target. This problem requires maintaining line-of-sight constraints while navigating to a moving target under thrust limitations. Figure \ref{fig:dockingres} presents the comparative trajectory analysis, control profiles, safety function evolution ($h_0$), and CLF behavior across all three control strategies and the combination of Stage 1 and 2. Again, initial states of failed trajectories where safety constraints cannot be satisfied are highlighted in red. The total control effort distributions are visualized in Figure \ref{fig:dockingcontrolthrust}, with corresponding statistical summaries provided in Table \ref{tableCC}.

The rendezvous results demonstrate consistent performance improvements similar to the cruise control case. Stage 1 achieves superior fuel efficiency compared to baseline ICCBF, reducing median fuel consumption from 0.72 kN to 0.64 kN (11.1\% reduction). The two-stage approach significantly expands mission feasibility, as within the 100 test cases, both ICCBF and Stage 1 alone fail in 12 instances (88\% success rate), while the combined Stage 1+2 framework reduces failures to only 5 cases (95\% success rate), recovering some cases that start with $x_0 \in \mathcal{S} \setminus \mathcal{C}^*$ representing an improvement in mission success rate when the initial state $\boldsymbol{x}_0 \in \mathcal{S}$. While the median fuel consumption remains the same going from Stage 1 to Stage 1+2, the combined approach exhibits extended upper tails in the fuel consumption distribution, which reflects the additional cost of recovering trajectories from the expanded feasible region $\mathcal{S} \setminus \mathcal{C}^*$ to $\mathcal{C}^*$.  Note that for failure cases during docking, the recommended approach is a safe abort and hold, followed by replanning.

\subsection{Comparative Analysis}
 Both test cases demonstrate improvements across key performance metrics such as fuel consumption and feasibility.  The cruise control problem exhibits slightly larger fuel savings, likely due to its lower-dimensional state space (2D vs 5D) and simpler constraint structure, allowing more effective RL optimization. 
 Furthermore, the rendezvous docking problem has the additional challenge of an original CBF $h_0$ that is independent of velocity, requiring the RL agent to learn velocity-dependency by itself to maintain constraint satisfaction. This, coupled with the higher-dimensional state space, makes the docking scenario significantly more challenging than cruise control. This is likely what causes the more modest fuel consumption reduction (11.1\% vs 12.5\%) observed in the rendezvous scenario compared to cruise control.

Th e introduction of Stage 2 reduces the failure rate of cruise control and docking, indicating that the two-stage framework can provide the opportunity for some initial states to retreat back into an inner safe set. Both cases show that with the two-stage RL framework, the expanded feasibility occurs at the cost of increased control effort variance. The extended upper tails in fuel consumption distributions ($P_{99}$ values of 0.68 and 1.00 for cruise control and rendezvous, respectively) represent the inherent cost of trajectory recovery from previously infeasible regions.

\section{Extended Application: 3D Uncooperative Target Inspection}\label{application}
Following the verification of the two-stage RL framework on the benchmarking test cases from \cite{agrawal2021safe}, it is applied to a three-dimensional inspection task adopted from \cite{doi:10.2514/1.G006126}, in which one spacecraft optimizes its trajectory to inspect an unknown non-cooperative resident space object (RSO). This section details the problem formulation, the framework implementation, and the obtained results.

\subsection{Problem Formulation}

\begin{figure}[hbt!]
    \centering    \includegraphics[width=\linewidth]{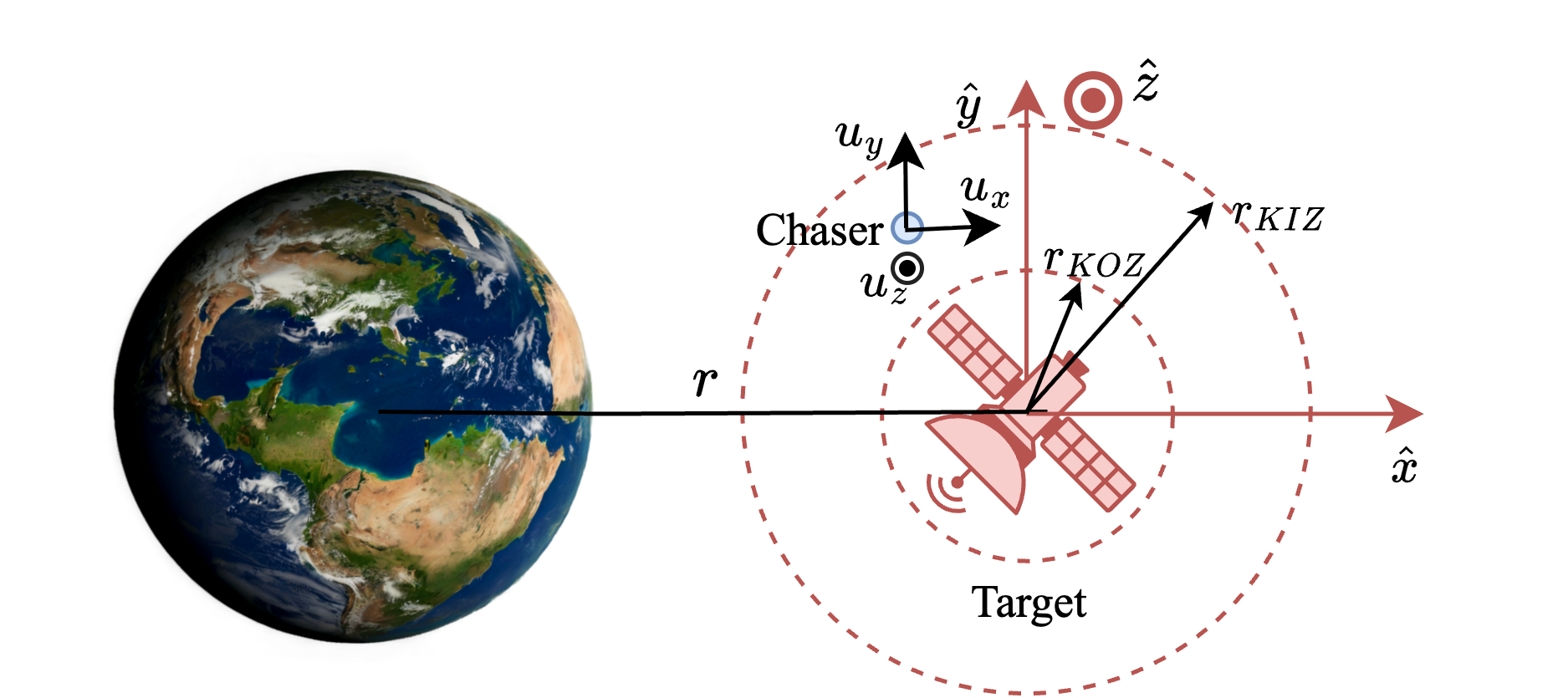}
    \caption{3D inspection problem}
    \label{fig:inspection}
\end{figure}

 In contrast to rendezvous and docking scenarios, the objective of this test case is not physical capture or contact but rather the maximisation of an inspection metric that accounts for observation quality.  The chaser spacecraft motion is modelled using the nonlinear relative motion equations in the LVLH frame, in 3 dimensions as 
\begin{equation}
\frac{d}{dt}
\begin{bmatrix}
p_x \\ p_y \\ p_z \\ v_x \\ v_y \\ v_z
\end{bmatrix}
=
\begin{bmatrix}
v_x \\
v_y \\
v_z \\
n^{2}p_x + 2n v_y + \dfrac{\mu}{r^{2}} - \dfrac{\mu(r+p_x)}{r_c^{3}} \\
n^{2}p_y - 2n v_x - \dfrac{\mu p_y}{r_c^{3}} \\
-\dfrac{\mu p_z}{r_c^{3}}
\end{bmatrix}
+ \frac{1}{m_c}
\begin{bmatrix}
0 \\ 0 \\ 0 \\ u_x \\ u_y \\ u_z
\end{bmatrix},
\label{eq:nonlinear-relative-3d}
\end{equation}

where $\mathbf{r}_c = [p_x, p_y, p_z]^T$ is the chaser’s relative position with respect to the RSO, $r_c = \sqrt{p_x^2 + p_y^2 + p_z^2}$ is the relative distance, $n = \sqrt{\mu/r^3}$ is the mean motion of the target orbit of radius $r = \SI{6771}{\kilo\meter}$, $m_c = \SI{50}{\kilo\gram}$ is the chaser mass, and $\mathbf{u} = [u_x, u_y, u_z]^T$ is the control acceleration vector, and $u_{max} =\SI{0.05}{\newton}$. The total time for inspection is set to be 48 h, and thrust and coast arcs are permitted during that time. The minimum and maximum $\Delta v$ for each of the burn arcs are $\SI{3}{\milli\meter\per\second}$ and  $\SI{90}{\milli\meter\per\second}$ and the maximum and minimum durations of the coast arcs are given as $t_{c} \in [t_{c_{\min}},  t_{c_{\max}}] = [\SI{1}{\hour}, \SI{3}{\hour}]$.  This test case is shown in detail in Figure \ref{fig:inspection}. 

The inspection performance is quantified by the cumulative score  
\begin{equation}\label{cgamma}
C_{\gamma} = \omega_{\gamma} \int_{t_0}^{t^{arc}} f(t)\,\frac{\pi - \gamma(t)}{\pi}\, dt,
\end{equation}
where $\gamma(t)$ denotes the instantaneous Sun–RSO–chaser angle, and $\omega_{\gamma}$ is a weighting parameter. The function $f(t)$ provides a distance-based weighting, defined as  
\begin{equation}
f(t) =
\begin{cases}
r_c/r_{\min}, & r_c < r_{\min}, \\
1, & r_{\min} \leq r_c \leq r_{\max}, \\
(r_{\max}/r_c)^3, & r_c > r_{\max},
\end{cases}
\end{equation}
with $r_{\min} = 50 \,\text{m}$ and $r_{\max} = 300 \,\text{m}$.   Note that this observability metric is from \cite{doi:10.2514/1.G006126}.

To ensure safety, two geometric constraints are enforced: (i) the chaser must remain within a keep-in zone (KIZ) of radius $r_{\text{KIZ}} = 1200 \,\text{m}$, and (ii) the trajectory must avoid penetration of a keep-out zone (KOZ) of radius $r_{\text{KOZ}} = 15 \,\text{m}$ surrounding the RSO.   As such, these constraints are expressed as $h_1(x)\geq 0$ and $h_2 (x) \geq 0$ where 
\begin{align}
    h_{{KOZ}_0}(x)  &= \frac{r_c^2 - r_{\text{KOZ}}^2}{r_{\text{KIZ}}^2 - r_{\text{KOZ}}^2}\\
    h_{{KIZ}_0}(x)  &= \frac{ r_{\text{KIZ}}^2 - r_c^2}{r_{\text{KIZ}}^2 - r_{\text{KOZ}}^2}
\end{align}
Note that $h_{KOZ}(x)$ and $h_{KIZ}(x)$ also have relative degree challenges as $L_g h(x) = 0$.

The inspection guidance problem can therefore be stated as follows: determine a feasible control policy $\mathbf{u}(t)$ that maximises the inspection score $C_{\gamma}$, subject to the nonlinear three-dimensional relative dynamics while maintaining the KIZ and KOZ constraints and adhering to the thrust constraints.

\subsection{Modified RL Framework}

The spacecraft inspection problem introduces additional complexity beyond the basic ICCBF framework described in Section \ref{methodology}, requiring RL to learn mission-specific parameters and control strategies. The inspection mission consists of alternating burn and coast phases, whose durations---denoted $t_b$ and $t_c$, respectively---significantly impact fuel efficiency and coverage performance. Thus, RL must learn to optimize these through actions $\eta_{b}, \eta_{c} \in [-1,1]$ where
\begin{align}
t_{b} &= t_{b_{\min}} + 0.5(\eta_{b} + 1)(t_{b_{\max}} - t_{b_{\min}}),
\label{eq:burn_duration} \\
t_{c} &= t_{c_{\min}} + 0.5(\eta_{c} + 1)(t_{c_{\max}} - t_{c_{\min}}).
\label{eq:coast_duration}
\end{align}
The burn duration bounds are determined by mission $\Delta v$ requirements: $t_{b_{\min}} = \frac{\Delta v_{\min}}{u_{\max}/m} = \frac{\SI{3}{\milli\meter\per\second}}{\SI{0.05}{\newton}/\SI{1000}{\kilo\gram}} = \SI{3}{\second}$ and $t_{b_{\max}} = \frac{\Delta v_{\max}}{u_{\max}/m} = \SI{90}{\second}$.

Additionally, the standard QP formulation that minimizes control effort to satisfy CBF constraints cannot directly optimize inspection metrics. To address this limitation, RL learns an additional inspection-oriented control component consisting of a unit direction vector $\hat{\boldsymbol{u}}_{RL} \in \mathbb{R}^2$
and a magnitude fraction $\lambda \in [0,1]$ representing the allocation of total control authority to inspection enhancement.  The inspection enhancement thrust is then calculated as:
\begin{equation}
\boldsymbol{u}_{RL} = \lambda \cdot u_{\max} \cdot \hat{\boldsymbol{u}}_{RL}.
\end{equation} This additional control is applied at the beginning of each RL step, allowing the agent to proactively influence the trajectory for improved inspection performance while the subsequent QP ensures safety.

In Stage 1,  QP formulation is now changed to 
\begin{equation}
    \begin{aligned}
&\underset{ \boldsymbol{u} \in \mathbb{R}^3, \; \delta}{\text{argmin}}
& & \frac{1}{2}  \boldsymbol{u}^T  \boldsymbol{u}+ p_2 \gamma_1 + p_3 \gamma_2 \\
&\text{subject to}
&& L_f  h_{1_{KOZ}} + L_g  h_{1_{KOZ}}\,  \boldsymbol{u} \geq -(\alpha_{RL,1} + \gamma_1)  h_{1_{KOZ}}, \\
&&& L_f h_{2_{KOZ}} + L_g h_{2_{KOZ}} \,  \boldsymbol{u} \geq -(\alpha_{RL,2} + \gamma_2) h_{2_{KOZ}}\\
&&& \parallel \boldsymbol{u}\parallel  \leq u_{\text{max}}
\end{aligned}
\label{QPinsp1}
\end{equation}
where
\begin{align}
    \label{hstage1ins}
    h_{1_{KOZ}}(\boldsymbol{x}) &=  b_{2,KOZ}(\boldsymbol{x}) \\
    h_{1_{KIZ}}(\boldsymbol{x}) &=  b_{2,KIZ}(\boldsymbol{x}) 
\end{align}
and $\boldsymbol{u} = \boldsymbol{u}_{RL} + \boldsymbol{u}_{QP}$. $\boldsymbol{u}_{QP}$ is the additional control added by the QP solution to maintain the CBFs.  $ b_{2,KOZ}(\boldsymbol{x})$ and $ b_{2,KIZ}(\boldsymbol{x})$ are the 2nd ICCBF layers for both the CBFs.  Here RL must determine $\alpha_{RL,1}$ and $\alpha_{RL,2}$ in addition to $\eta_c$, $\eta_b$ and $\boldsymbol{u}_{RL}$.

In Stage 2, $h_{RL}$ is learnt and used to define a new CBF 
\begin{equation}\label{hstage2ins}
    h (\boldsymbol{x}) = \frac{1}{2}(h_{{KOZ}_0} + h_{{KIZ}_0}) + \overline{h}_0 h_{RL}(\boldsymbol{x}).
\end{equation}
 Note that $\overline{h}_0 =   \frac{1}{2}(\overline{h}_{{KOZ}_0} + \overline{h}_{{KIZ}_0})$. In Stage, 2 RL must determine $h_{RL},  \alpha_{RL,1}$ and $\alpha_{RL,2}$ in addition to $\eta_c$, $\eta_b$ and $\boldsymbol{u}_{RL}$. The Stage 2 QP is formulated as 
\begin{equation}
    \begin{aligned}
&\underset{ \boldsymbol{u} \in \mathbb{R}^3, \; \delta}{\text{argmin}}
& & \frac{1}{2}  \boldsymbol{u}^T  \boldsymbol{u} + p \gamma \\
&\text{subject to}
&& L_f  h + L_g  h\,  \boldsymbol{u} \geq -(\alpha_{RL} + \gamma)  h, \\
&&& \parallel  \boldsymbol{u}\parallel  \leq  u_{\text{max}}
\end{aligned}
\label{QPinsp2}
\end{equation}
For both stages, the reward is calculated as 
\begin{equation}\label{inspectionreward}
\begin{aligned}
R_j ={}& - c_{h_1}\,\text{max}(0,-h_{{KOZ}_0}(\boldsymbol{x}_j)) \\
       & - c_{h_2}\,\text{max}(0,-h_{{KIZ}_0}(\boldsymbol{x}_j)) - c_i C_{\gamma}\, .
\end{aligned}
\end{equation}
to include both CBFs. Note that $ c_{h_1} = 1.0$, $ c_{h_2} = 1.0$ and $c_i = 0.1$ are user set reward gains and $C_{\gamma}$ is calculated via Eq.~\eqref{cgamma}. The extended RL environment for the inspection problem is given in Algorithm \ref{inspection}, with the network architecture and hyperparameters for this case given in Table \ref{insphyp}.  These parameters and network architecture were chosen empirically by running small training batches with varying settings in a simple parameter sweep.

\begin{table}[hbt!]
\centering 
\caption{RL hyperparameters for the inspection problem}
\renewcommand{\arraystretch}{1.3}
\label{insphyp}
\begin{tabular}{lcc}
\hline
Parameter                   & \multicolumn{2}{c}{Inspection}                              \\ \cline{2-3} 
                            & \multicolumn{1}{c|}{Stage 1}            & Stage 2            \\ \hline
Learning Rate               & \multicolumn{1}{c|}{$10^{-4}$ (const.)} & $10^{-4}$ (decay.) \\
Batch Size                  & \multicolumn{1}{c|}{64}                 & 256                \\
Rollout Steps               & \multicolumn{1}{c|}{1,280}              & 1,280              \\
Training Epochs             & \multicolumn{1}{c|}{10}                 & 10                 \\
Discount Factor ($\gamma$)  & \multicolumn{1}{c|}{0.95}               & 0.95               \\
GAE Lambda ($\lambda$)      & \multicolumn{1}{c|}{0.99}               & 0.99               \\
Clip Range                  & \multicolumn{1}{c|}{0.2}                & 0.2                \\
Entropy Coefficient         & \multicolumn{1}{c|}{0.01}               & 0.01               \\
Initial Std. Dev.           & \multicolumn{1}{c|}{0.2}                & 0.21                \\
SDE & \multicolumn{1}{c|}{Yes}                & No                 \\
Network Layers              & \multicolumn{1}{c|}{4}                  & 4                  \\
Neurons per Layer           & \multicolumn{1}{c|}{64}                 & 64                 \\
Activation Function         & \multicolumn{1}{c|}{Tanh}               & Tanh        \\ \hline       
\end{tabular}
\end{table}

\begin{algorithm}
\caption{RL environment for the inspection problem}
\label{inspection}
\textbf{Input:}
Actions $\mathbf{a}_j=[\alpha_{1,RL,j},\alpha_{2,RL,j}, \eta_{b,j},\eta_{c,j},\hat{\boldsymbol{u}}_{RL,j}, \lambda_j]$ (Stage 1) or $[h_{\mathrm{RL},j},\alpha_{RL,j},\eta_{b,j},\eta_{c,j},\hat{\boldsymbol{u}}_{RL,j}, \lambda_j]$ (Stage 2); state $\boldsymbol{x}_j$; current time $t$; final time $t_f$; .

\begin{algorithmic}[1]
\State Compute burn time $t_b$ and thrust time $t_c$ using Eq. \eqref{eq:burn_duration} and Eq. \eqref{eq:coast_duration}.
\State Compute $h(\boldsymbol{x}_j)$ using Eq.~\eqref{hstage1ins} for Stage 1 or Eq.~\eqref{hstage2ins} for Stage 2.
\State Compute $\partial h/\partial \boldsymbol{x}$ and the corresponding Lie derivatives.
\State Compute the ZOH control $\boldsymbol{u}_j$ by solving the QCQP in Eq.~\eqref{QPinsp1} for Stage 1 or Eq.~\eqref{QPinsp2} for Stage 2.  
\State Forward propagate the system from $t$ to $t+ t_b$ under $\boldsymbol{u}_j +\boldsymbol{u}_{RL}$, followed by a propagation from $t+t_b$ to $t + t_b + t_c$ under ballistic conditions using Eq.~\eqref{eq:nonlinear-relative-3d} to obtain $\boldsymbol{x}_{j+1}$.
\State Update $t \gets t+t_c + t_b$.
\State Compute the reward $R_{j+1}$ using~\eqref{inspectionreward}.
\State Calculate the propagated actor input state:
\[
\boldsymbol{S}_{j+1} =
\begin{cases}
[\boldsymbol{x}_{j+1}], 
& \text{Stage 1}, \\[6pt]
\begin{aligned}
&[\boldsymbol{x}_{j+1},\; L_g h(\boldsymbol{x}_{j+1}),\; L_f h(\boldsymbol{x}_{j+1}), \\
&\quad h_0(\boldsymbol{x}_{j+1}),\; V(\boldsymbol{x}_{j+1})],
\end{aligned}
& \text{Stage 2}.
\end{cases}
\]
\State Scale $\boldsymbol{S}_{j+1}$ to obtain $\boldsymbol{S}_{j+1}^\ast$ using~\eqref{actorin}.
\State \textbf{Output:} $\boldsymbol{S}_{j+1}^\ast,\ R_{j+1}$.
\If{$t = t_f$}
  \State \textbf{Stop.}
\Else
  \State \textbf{Update:} $j \gets j+1$ and repeat from Step~1.
\EndIf
\end{algorithmic}
\end{algorithm}

\subsection{Inspection Results}
The inspection mission performance is evaluated across a range of initial conditions that represent realistic approach scenarios. Initial states are selected such that $\boldsymbol{x}_0 \in \mathcal{S}$, with the initial spacecraft positioned placed at varying distances from the target. The initial conditions are parameterized as
$\boldsymbol{x}_0 = \begin{bmatrix}
r_0 & 0 & 0 & 0 & -2nr_0 & 0.431
\end{bmatrix}^T$, where $r_0 \in [50, 300]$ m represents the initial radial distance from the target.

\begin{figure*}[hbt!]
    \centering
    \includegraphics[width=\linewidth]{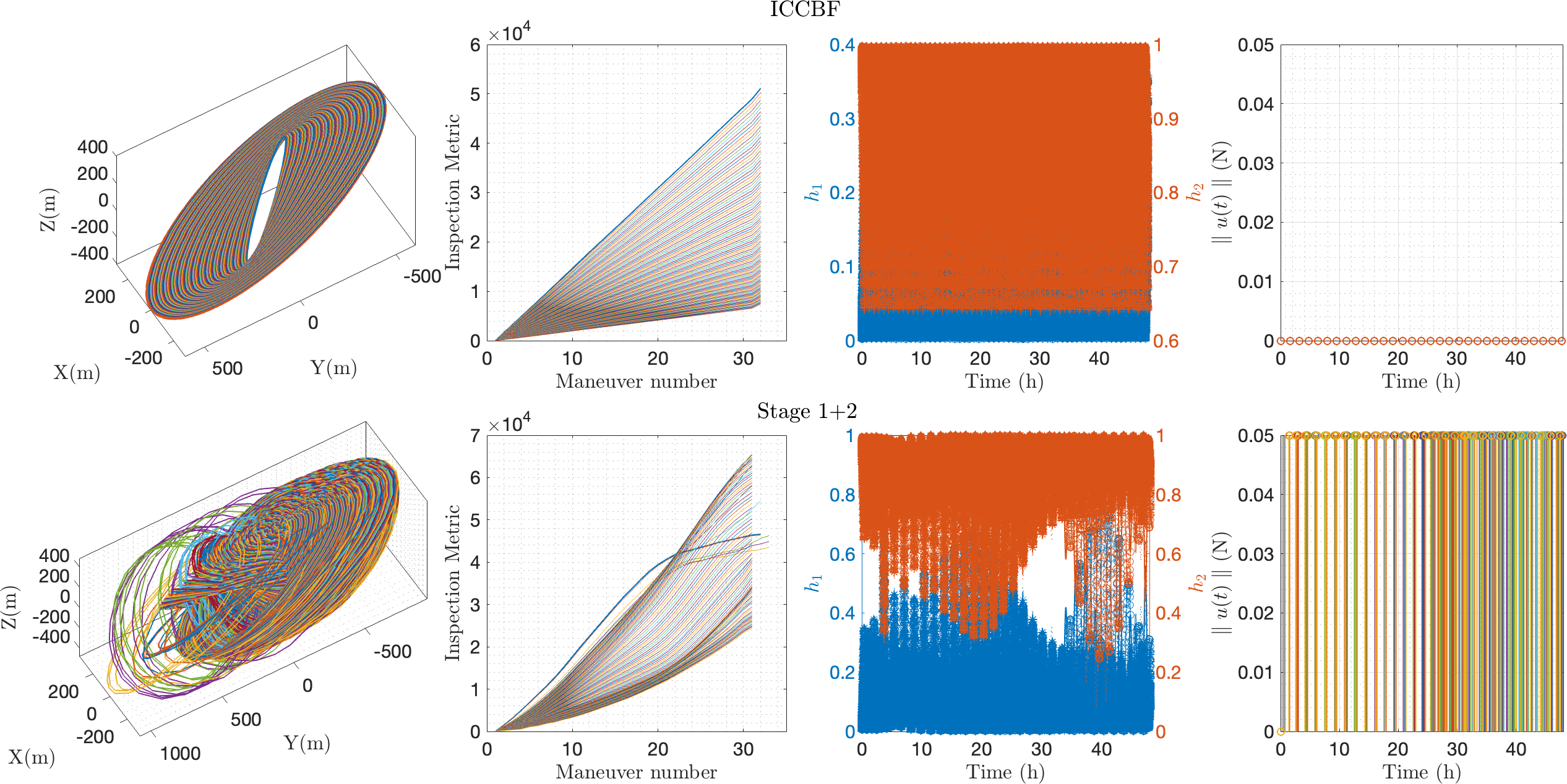}
    \caption{3D inspection Results. Comparative analysis of trajectory performance between ICCBF and Stage 1-2 combination, inspection metric evolutions, safety constraints ($h_1$ and $h_2$), and control profile across time.}
    \label{fig:insp}
\end{figure*}

\begin{figure}[hbt!]
    \centering
    \includegraphics[width=0.5\linewidth]{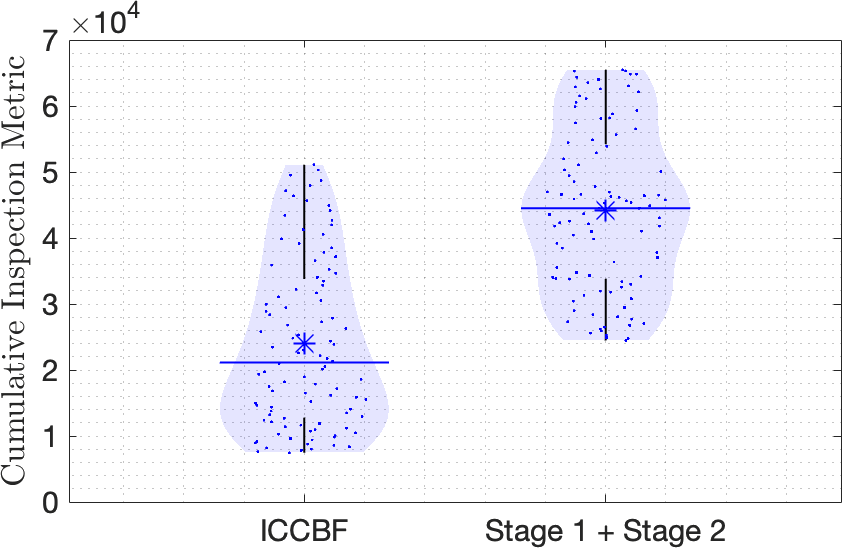}
    \caption{Cumulative inspection metric distribution for the inspection problem}
    \label{fig:metricinspection}
\end{figure}

A batch of 100 MC simulations are conducted under each control strategy, and the the baseline is provided by the ICCBF implementation with fixed timing parameters ($\eta_b = \eta_c = 0$, corresponding to nominal burn and coast durations) and no inspection enhancement ($\boldsymbol{u}_{\text{insp}} = \boldsymbol{0}$).  The inspection performance results are presented in Figure \ref{fig:insp}, which shows representative trajectories, thrust usage, and safety constraint adherence for ICCBF and combined Stage 1 and 2 approaches. Figure \ref{fig:metricinspection} displays the cumulative distribution of inspection metrics, demonstrating the improvement achieved by the RL-enhanced approach across the test population. The quantitative inspection metric comparison is summarized in Table \ref{tableinsp}, which presents statistical measures of the total inspection metric score accumulated over the 48-h mission duration.

\begin{table}[hbt!]
\renewcommand{\arraystretch}{1.3}
\centering
\begin{tabular}{lcc} \hline
{CBF Strategy} & \textbf{$\mu$ $\pm$ $\sigma$} [$\times 10^4$] & [Q$_1$, Q$_2$, Q$_3$, P$_{99}$] [$\times 10^4$] \\ \hline 
ICCBF & 2.39 $\pm$ 1.27 & [1.28, 2.11, 3.37, 5.07] \\
Stage 1+2& 4.42 $\pm$ 1.23 & [3.38, 4.45, 5.42, 6.54] \\ \hline
\end{tabular}
\caption{Statistical summary of total inspection metric.}
\label{tableinsp}
\end{table}

The results demonstrate significant performance improvements with the RL-enhanced methodology. The average inspection metric score increases by 85\%, from 2.39$\times 10^4$ to 4.42$\times 10^4$, representing nearly doubled inspection capability over the 48-hour mission. All quartile measures show substantial increases in the combined Stage 1+2 case.

\subsection{Performance Analysis}
With the current initial stage distribution being well within the safety zone, Stage 2 (RL-residual CBF) was required for trajectory recovery in only 3 out of the 100 simulations. While this illustrates that Stage 1 is capable of handling nearly all trajectories if the initial state is well within the safety zone, it limits the testing of the recovery ability of Stage 2. Thus, future work should include testing with initial conditions near or on the KIZ boundary to better evaluate Stage 2 recovery capabilities.

While the initial inspection results presented here show promise, several avenues for improvement remain. Exploration of deeper network architectures represents a crucial component of future work, as it may yield better performance for this complex multi-objective scenario. Additional research directions include extended training with more challenging initial conditions, comprehensive analysis of trade-offs between safety margins and inspection performance in boundary regions, and investigation of alternative reward structures to better balance competing objectives.

\section{Conclusion}\label{conclusions}

This work presents a unified two-stage RL framework that addresses key limitations of ICCBF for safety-critical spacecraft control. Stage 1 learns adaptive class-$\mathcal{K}_\infty$ parameters that reduce ICCBF conservatism while maintaining safety guarantees, while Stage 2 expands the feasible operating region of ICCBF by learning a residual barrier function for previously abandoned states. The methodology demonstrates significant improvements across three problems: 12.5\% fuel reduction in cruise control, 11.1\% fuel savings in spacecraft docking, and 85\% inspection metric increase in 3D spacecraft RSO inspection tasks. The framework preserves the computational efficiency and theoretical safety guarantees of the CLF-CBF-QP structure while embedding long-term trajectory optimization awareness, showing that RL can effectively make CBFs non-greedy without sacrificing real-time performance or safety certification.




\bibliographystyle{abbrv}        
\bibliography{refs}        

\end{document}